\documentclass[12pt]{article}
\usepackage{ amsfonts,
amsmath, amssymb,amsgen, amsthm, amscd}
\input psfig.sty

\def\Diff{\mathop{\rm Diff}\nolimits}
\def\A{\mathop{\rm A}\nolimits}
\def\BT{\mathop{\rm BT}\nolimits}

\def\GL{\mathop{\rm GL}\nolimits}

\def\GV{\mathop{\rm GV}\nolimits}

\def\Hopf{\mathop{\rm Hopf}\nolimits}
\def\Id{\mathop{\rm Id}\nolimits}
\def\Im{\mathop{\rm Im}\nolimits}
\def\Re{\mathop{\rm Re}\nolimits}
\def\det{\mathop{\rm det}\nolimits}
\def\ad{\mathop{\rm ad}\nolimits}
 
\def\Ker{\mathop{\rm Ker}\nolimits}

\def\log{\mathop{\rm log}\nolimits}

\def\Tr{\mathop{\rm Tr}\nolimits}
\def\SL{\mathop{\rm SL}\nolimits}
\def\PSL{\mathop{\rm PSL}\nolimits}

\def\Ac{{\cal A}}

\def\Ec{{\cal E}}
\def\Fc{{\cal F}}
\def\Gc{{\cal G}}
\def\Hc{{\cal H}}

\def\Lc{{\cal L}}
\def\Mc{{\cal M}}

\def\Rc{{\cal R}}
\def\Sc{{\cal S}}

\def\a{\alpha}
\def\b{\beta}
\def\d{\delta}
\def\D{\Delta}
\def\g{\gamma}
\def\G{\Gamma}
\def\lb{\lambda}
\def\Lba{\Lambda}
\def\om{\omega}
\def\Om{\Omega}
\def\s{\sigma}
\def\Si{\Sigma}
\def\t{\theta}
\def\ve{\varepsilon}
\def\vp{\varphi}

\def\ab{\bf a}

\def\mb{\bf m}
\def\xb{\bf x}
\def\yb{\bf y}
\def\zb{\bf z}
\def\Bb{\bf B}
\def\Eb{\bf E}
\def\Sb{\bf S}
\def\0b{\bf 0}

\def\bash{\backslash}
\def\bu{\bullet}
\def\fl{\forall}
\def\ify{\infty}
\def\mpo{\mapsto}
\def\nb{\nabla}
\def\ot{\otimes}

\def\ra{\rightarrow}

\def\sbs{\subset}

\def\ts{\times}
\def\wdg{\wedge}
\def\wt{\widetilde}

\newcommand{\FH}{H}

\newcommand{\Fa}{\mathfrak{a}}

\newtheorem{theorem}{Theorem}

\newtheorem{proposition}[theorem]{Proposition}
\newtheorem{lemma}[theorem]{Lemma}
\newtheorem{corollary}[theorem]{Corollary}

\def\limproj{\mathop{\oalign{lim\cr
\hidewidth$\longleftarrow$\hidewidth\cr}}}
\def\limind{\mathop{\oalign{lim\cr
\hidewidth$\longrightarrow$\hidewidth\cr}}}

\def\build#1_#2^#3{\mathrel{
\mathop{\kern 0pt#1}\limits_{#2}^{#3}}}

\numberwithin{equation}{section}
\begin{document}

\title{{\bf  Modular Hecke Algebras and their Hopf
Symmetry}}
\author{Alain Connes \\
        Coll\`ege de France \\
        3 rue d'Ulm \\
        75005 Paris, France \\
\and
        Henri Moscovici\thanks{Research
    supported by the National Science Foundation
    award no. DMS-9988487.} \\
    Department of Mathematics \\
    The Ohio State University \\
    Columbus, OH 43210, USA
    }

\date{ \ }

\maketitle

\centerline{\it Dedicated to Pierre Cartier} 

\begin{abstract}

We introduce and begin to analyse a class of algebras, 
associated to congruence
subgroups, that extend both 
the algebra of modular forms of all levels and the ring of classical 
Hecke operators. At the intuitive level,
these are algebras of `polynomial coordinates'
for the `transverse space' of lattices
modulo the action of the Hecke correspondences. Their 
underlying symmetry is shown to
be encoded by the same Hopf algebra that controls the transverse
geometry of codimension 1 foliations. Its action is shown
to span the `holomorphic tangent space' of the noncommutative space,
and each of its three basic Hopf cyclic cocycles acquires a specific
meaning.
The Schwarzian 1-cocycle gives an inner derivation implemented by the 
level 1 Eisenstein series of weight 4. The Hopf cyclic 2-cocycle 
representing the transverse fundamental class provides a natural 
extension of the first Rankin-Cohen bracket to the modular Hecke 
algebras. Lastly, the Hopf cyclic version of the Godbillon-Vey cocycle 
gives rise to a 1-cocycle on $\PSL(2, \mathbb{Q})$ with values in Eisenstein series 
of weight 2, which, when coupled with the `period' cocycle, yields a 
representative of the Euler class.
  
\end{abstract}

\section*{Introduction}

The aim of this paper is to introduce and begin to analyse
a  class of algebras, called modular Hecke algebras, which 
encode the two \textit{a priori} unrelated structures 
on modular forms, namely the algebra structure given by the pointwise product
 on the one hand,
and the action of the Hecke operators on the other.
We associate to any congruence subgroup $\G$
of $\, \SL(2, \mathbb{Z})$ a crossed product algebra $\Ac(\G)$,
the \textit{modular Hecke algebra} of level $\G$,  which is a
direct extension of both the ring of classical Hecke operators and
of the algebra $\Mc(\G)$ of $\G$-modular forms. With $\Mc$
denoting the algebra of modular forms of arbitrary level, the
elements of $\Ac(\G)$ are maps with finite support
$$
 F: \G\bash \GL^+ (2, \mathbb{Q}) \ra \Mc \, ,
\qquad \G  \a \mapsto F_{\a} \in \Mc  \, ,
$$
satisfying the covariance condition
\begin{equation}
F_{\a \g} \, = \, F_{\a} \vert \g \, , \qquad \fl \a \in \GL^+ (2,
\mathbb{Q}) \, , \g \in \G \, \nonumber
\end{equation}
and their product is given by convolution. In the simplest case
$\G(1)=\,\SL(2, \mathbb{Z})$, the elements of $\Ac(\G(1))$ are
encoded by a finite number of modular forms $f_N \in
\Mc(\G_0(N))\,$ of arbitrary high level and the product operation
is non-trivial. The algebra $\Ac(\G)$ acts on $\Mc(\G)$ extending
the action of classical Hecke operators, and the `cuspidal'
elements form a two-sided ideal of $\Ac(\G)$.\\

\noindent Our starting point is the basic observation that the
Hopf algebra $\Hc_1$, which was discovered (cf.~\cite{CM1}) in the
analysis of codimension $1$ foliations as underlying symmetry of the
`transverse' geometry, admits a natural action on the crossed
product algebras  $\Ac(\G)$. This
action is sufficiently non-trivial to span the `tangent space'
of the non-commutative space $Q(\G)$ with coordinate ring $\Ac(\G)$,
and is thus a key ingredient in the understanding of the 
geometry of $Q(\G)$. At an intuitive level these non-commutative spaces 
are obtained from the quotient of the space of lattices by
the Hecke correspondences, while  $\Ac(\G)$ is the
simplest  algebra
of polynomial coordinates on $Q(\G)$.\\

\noindent To describe the action of $\Hc_1$,
we recall that, as an algebra, $\Hc_1$
coincides with the universal enveloping algebra of the Lie algebra
with basis $\{ X,Y,\d_n \, ; n \geq 1 \}$ and brackets
\begin{equation*}
[Y,X] = X \, , \, [Y , \d_n ] = n \, \d_n \, , \, [X,\d_n] =
\d_{n+1} \, , \, [\d_k , \d_{\ell}] = 0 \, , \quad n , k , \ell
\geq 1 \, ,
\end{equation*}
while the coproduct which confers it the Hopf algebra structure is
determined by the identities
\begin{eqnarray}
\D \,  Y &=& Y \ot 1 + 1 \ot Y \, , \quad \D \,  \d_1 = \d_1 \ot 1
+ 1 \ot \d_1 \, , \nonumber
 \\
\D \,  X &=& X \ot 1 + 1 \ot X + \d_1 \ot Y \, , \nonumber
\end{eqnarray}
together with the property that $\, \D : \Hc_1 \ra \Hc_1 \ot \Hc_1
\, $ is an algebra homomorphism. The action of $X$ on $\Ac (\G)$
is given by a classical operator going back to Ramanujan
(\cite{Ra}, \cite{Se}), which corrects the usual differentiation
by the logarithmic derivative of the Dedekind eta function $ \eta
(z)$. The action of $Y$ is given by the standard grading by the
weight (the Euler operator) on modular forms. Finally, $\d_1$ and
its higher `derivatives' $\d_n$ act by generalized cocycles on $
\GL^{+} (2 , \mathbb{Q})$ with values in modular forms.
\medskip

\noindent The main result of the present paper is a rather
complete understanding of the above action of $\Hc_1$, viewed as
the underlying symmetry of the noncommutative space $Q(\G)$
 obtained from the quotient of the space of lattices by
the Hecke correspondences. The
picture that emerges is that of a surprisingly close analogy
between the action of Hecke operators on modular forms and the
action of a discrete subgroup of $\Diff (S^1)$ on polynomial
functions on the frame bundle of $S^1$. The role of the angular
variable $\theta \in \mathbb{R}/2 \, \pi \,\mathbb{Z}$ is assumed
by the simplest Eichler integral, namely the primitive $Z \in
\mathbb{C}/\Lambda$ of the holomorphic differential form $\frac{2 \pi i}{6}\,\eta^4
\, dz$ on the
 elliptic curve $X'(1) = \Gamma'(1)\bash {\FH}^* $,
 where $\Gamma'(1)$ is the commutator subgroup of
 $\Gamma(1)$.
The lattice $2 \, \pi \,\mathbb{Z}\subset \mathbb{R}$ is replaced
by the equilateral lattice $\Lambda \subset \mathbb{C}$ of periods
of $Z$. The part of the circle
$$ e^{i \theta} =\,
\cos \theta \, + \, i \, \sin \theta \, , \qquad x^{2} \,+ \,
y^{2} = \, 1
$$
is thus played by the genus $1$ curve $X'(1) \sim
\mathbb{C}/\Lambda$ with Weierstrass parametrization and
equation
$$
\left( \wp_{\Lambda}(Z),\, \wp'_{\Lambda}(Z)\right)
\,=\,\left(\sqrt[3]{j},\, \frac{-2 E_6}{\eta^{12}}\right) \, , \qquad
y^2  = \, 4 \,(x^3 \, - \,1728) \, 
$$
using traditional notations. The diffeomorphisms $\,
\phi(\theta)\, $ are replaced by Hecke transformations
$\, Z|_0\, \gamma \,$, while the Jacobian $\,
\phi'(\theta) \,$ becomes $\, J(\gamma) =\,\frac{d
Z|_0\, \gamma}{d Z} $. Accordingly, the higher cocycles
$ \,\left(\frac{d}{d\theta}\right)^n \, \log
\phi'(\theta) \,$ which allow us to define the action of $\d_n$ take
the form $\,\left(\frac{d}{dZ}\right)^n \, \log
J(\gamma) $. Furthermore, the usual projective structure on
$S^{1}$, corresponding to $\, t \, = \, \tan
(\frac{\theta}{2}) \,$, which in terms of $\theta$ is given by the
quadratic differential
\begin{equation*}
    \rho \, : = \, \{t\, ; \theta \} \, d \theta^2 \, = \,
    \frac{1}{2} \, d \theta^2 \, ,
\end{equation*}
has as its modular counterpart the rational projective structure
on $X'(1)$ given by the quadratic differential
\begin{equation*}
  \varpi' \, = \, \frac{x \,  dx^2 }{2 \, y^2} \, = \,
 \frac{x \,  dx^2 }{8 \,(x^3-1728)} \, .
\end{equation*}
This Schwarzian relation between $Z$ and
$\varpi'$ is similar to the classical relation between the
modular invariant $j$ and
$E_4$ in the case of genus zero \cite{Mc}.
\medskip

\noindent With the above dictionary in mind, one can proceed to
transfer transverse geometry concepts and constructions to the
setting of modular forms (cf. Theorem \ref{trans} below). 
The basic tool is the cyclic cohomology of Hopf
algebras developped in (\cite{CM1}). For $\Hc_1$
one gets three basic cyclic
classes, $\delta'_2$, $\delta_1$ and $F$ which 
in the
original action of $\Hc_1$ on the crossed product of functions on
the frame bundle on $S^1$ by a discrete subgroup of $\Diff (S^1)$
correspond respectively to 
\[
\begin{array}{cc}
\hbox{ Schwarzian} &\hbox{$\delta'_2:= \delta_2 - \frac{1}{2}  \delta_1^2 $} \
 \\
 & \\
\hbox{ Godbillon-Vey Class} &\hbox{$\delta_1$} \\
 & \\
\hbox{Transverse Fundamental Class} &\hbox{$F := X \ot Y - Y \ot X - \d_1 \,  Y \ot Y \,$ .
} \\
 \end{array}
\]
In the present paper we shall compute the  above 
three classes in the action of  $\Hc_1$ on the modular
Hecke algebras.\\

 \noindent We first show
that, under the present  action of $\Hc_1$, the element
$\delta'_2$ is represented by an inner derivation implemented by
the above mentioned Eisenstein series of weight 4 and level 1,
denoted here by $\om_4$. We then prove that no inner perturbation of
the action by a 1-cocycle can annihilate the Schwarzian
$\tilde{\delta'_2}$, and that the freedom one has in modifying the
action of the Hopf algebra $\Hc_1$ by a $1$-cocycle exactly
changes the restriction of the action of $X$ on modular forms and
the value of $\om_4$ as in the data used by Zagier~\cite{Za} to
define `canonical' Rankin-Cohen algebras, where moreover
the same $\om_4$  appears with a slightly
different normalization.\\

 \noindent We then show that the image
under the canonical map from the Hopf cyclic cohomology of
$\Hc_{1}$ to the Hochschild cohomology of the algebras $\Ac(\G)$
of the Hopf $2$-cocycle encoding the `transverse fundamental
class' gives the natural extension of the first Rankin-Cohen
bracket for modular forms to $\Ac(\G)$. Actually, in a sequel to
this paper~\cite{CM4}, we show that all the Rankin-Cohen brackets
have natural extensions to modular Hecke algebras, which moreover
are determined by a universal deformation formula 
based on $\Hc_1$. The
mere existence of the modular Hecke algebras is sufficient to
prove in full generality the universal associativity property.
\medskip

\noindent  While the transverse fundamental class generates the
even component of the periodic Hopf cyclic cohomology of $\Hc_1$,
the odd component is generated by the Hopf cyclic version of the
Godbillon-Vey class, namely $\delta_1$. In the context of the present action of
$\Hc_1$, this class gives rise to a `transverse' $1$-cocycle on
$\GL^+ (2, \mathbb{Q})$ with values in Eisenstein series of weight
$2$, measuring the lack of $\GL^+ (2, \mathbb{Q})$-invariance of
the connection $\nabla$ corresponding to $X$. When coupled with
the obvious extension of the classical period cocycle, it yields
an analogue of the Bott-Thurston group $2$-cocycle~\cite{Bo}, that
represents the Euler class $\,e \in \, H_{\rm d}^{2} \, (\SL (2,
\mathbb{R}), \, \mathbb{R}) \,$ (cf. Theorem \ref{m1}). By a
theorem of Borel and Yang \cite{BY}, the restriction of this class
to $\, H^{2} \, (\SL (2, \mathbb{Q})^{\d}, \, \mathbb{R}) \,$ is
nontrivial and generates $\, H^{2} \, (\SL (2, \mathbb{Q})^{\d},
\, \mathbb{R}) \,$.
\medskip

\noindent 
Our last result (Theorem \ref{m2}) gives an arithmetic
representation for the rational Euler class $\,e \in \, H^{2} \,
(\SL (2, \mathbb{Q})^{\d} , \, \mathbb{Q}) \,$ in terms of
generalized Dedekind sums \cite{De}, \cite{Me}.  
 The original Dedekind  sums already  arose 
 in several different topological settings,
as amply illustrated  by Atiyah~\cite{A} and
also by Kirby and Melvin~\cite{KM}. 
\medskip

\noindent In Section 1 we introduce the modular Hecke algebra
$\Ac(\G)$ for any congruence subgroup $\G \sbs \SL (2,
\mathbb{Z})$. Its elements are form-valued Hecke operators; the
classical Hecke algebra is the subalgebra of elements of degree
$0$ for the natural grading given by the weight.

\noindent Section 2 recalls the Hopf algebra $\Hc_1$, in its
original context of codimension $1$ transverse geometry, and
describes its most significant Hopf cyclic cocycles, that
correspond to the Godbillon-Vey class, to the Schwarzian
derivative, and to the transverse fundamental class.

\noindent In Section 3 we define and analyze in great detail the
Hopf action of $\Hc_1$ on the crossed product of modular forms by
Hecke operators. We show that the  Schwarzian  $1$-cocycle gives
an inner derivation. We relate the conjugates by inner
perturbations of the Hopf action of $\Hc_1$ to the data used by
Zagier~\cite{Za} to define canonical Rankin-Cohen algebras. We
also show that the  Hopf cyclic $2$-cocycle representing the
transverse fundamental class  provides a natural extension of the
first Rankin-Cohen bracket to the crossed product algebras. The
section concludes with a precise formulation of the geometric
analogy with the case of $\Diff(S^1)$. In particular, we give a
simple geometric explanation for the above formula for the
quadratic differential $\varpi$.

\noindent In Section 4 we  introduce the `transverse'  $1$-cocycle
on $\GL^+ (2, \mathbb{Q})$ with values in Eisenstein series of
weight $2$, as a partial analogue of the Godbillon-Vey cocycle. We
then complete the construction of its full analogue using the
`period' cocycle and prove that it represents the Euler class of
$\, \SL (2, \mathbb{R})$. Finally, we conclude by expressing the
Euler class of $\, \SL (2, \mathbb{Q})$ in terms of generalized
Dedekind sums.

\noindent Appendix A recalls the definition of cyclic cohomology
for Hopf algebras, while Appendix B provides the details for the
interpretation of the Godbillon-Vey class as a Hopf cyclic
cohomological class.
\bigskip

\noindent We wish to express our gratitude to Don Zagier --- this
paper has greatly benefited from his  course on modular forms at
Coll\`ege de France \cite{Zagier} as well as from his comments.
\bigskip

\tableofcontents

\section{Modular Hecke Algebras}

Modular forms, regarded as lattice functions, may be viewed as
natural coordinates on the lattice line bundle over the moduli
space of isomorphism classes of elliptic curves over $\mathbb{C}$.
We recall that every elliptic curve is isomorphic -- as a complex
manifold -- to a 1-dimensional complex torus $\, E = \mathbb{C} /
\Lba \,$, where $\, \Lba \,$ is a lattice in $\mathbb{C}$\, ;
furthermore, two lattices $\, \Lba $ and $ \Lba' \,$ define
isomorphic curves $\, E \simeq E' \,$ iff they are homothetic,
$$ \Lba' \, = \, \lb \Lba \, , \qquad \text{for some} \quad \lb \in
\mathbb{C}^{\ts} \, .
$$
Since one can always represent a lattice $\, \Lba \,$ in the form
$$
 \Lba = \mathbb{Z}\om_1 + \mathbb{Z} \, \om_2 \, , \qquad
  \Im \, \frac{\om_1}{\om_2} > 0 \, ,
$$
up to the obvious action of the group $\, \G(1) = \SL (2 ,
\mathbb{Z}) \,$, it follows that the isomorphism classes of
elliptic curves are parametrized by the quotient $\, \G(1) \bash H
\,$, where $\, {\FH} = \{ z \in \mathbb{C} , \Im (z) > 0 \}$ is
the upper half-plane. Thus, the set $\Lc$ of all lattices in
$\mathbb{C}$ defines a line bundle
$$
   \mathbb{C}^{\ts} \ra \Lc \ra \G(1) \bash {\FH} \, ,
$$
and the line bundle $ \Lc^{-2}$ is canonically isomorphic to the
(complex) $1$-jet bundle of $\G(1) \bash {\FH} $.
\medskip

\noindent Any lattice function which is homogeneous of weight
$2k$,
$$
F (\lb \Lba) = \lb^{-2k} F(\Lba) \, , \qquad \lb \in
\mathbb{C}^{\ts} \, ,
$$
is automatically of the form
$$
 F (\mathbb{Z}\om_1 + \mathbb{Z} \om_{2}) =
\om_2^{-2k} \, f \left( \frac{\om_1}{\om_2} \right)
$$
with the function $f$ satisfying the modularity property
\begin{equation}
f|_{2k} \, \g  \, = f \, , \qquad \fl \, \g \in \G(1) \, = \, \SL
(2 , \mathbb{Z}) \, ,
\end{equation}
where we used the standard `slash operator' notation for the
action of
$$
G^{+} (\mathbb{R}) := \GL^{+} (2, \mathbb{R})
$$
on functions on the upper half plane:
\begin{eqnarray} \label{act1}
f|_{k}\, \a \, (z) \, &=& \, \det (\a)^{k/2} \,
f (\a \cdot z) \, j(\a , z)^{-k} \, , \\
\a = \begin{pmatrix}a  &b\\
c   &d  \end{pmatrix} \in G^{+} (\mathbb{R}) , \quad
 \a \cdot z &=& \frac{az+b}{cz+d}
 \quad \hbox{and} \quad
j (\a, z) \, = \, cz+d \, . \nonumber
\end{eqnarray}
In particular, $\, f (z + 1) = f(z) , \, z \in \mathbb{C} \,$,
hence it induces a function on the punctured unit disc $\,
f_{\ify}(q)  \, , \quad q = e^{2  \pi i z} \, $. Such a function
$f \,$ is called a {\it modular form of weight} $2k$ if it is
holomorphic on ${\FH}$ and also at $\ify$, that is $f_{\ify}$ is

holomorphic at $ \, q = 0 \,$\, ; $\, f \,$is called a {\it cusp
form} if in addition $ \, f_{\ify}(0) = 0 \,$. We denote by
$$
\Mc (\G(1)) := \Si^{\oplus} \, \Mc_{2k} (\G(1)) \, , \quad
\hbox{resp.} \quad \Mc^{0} (\G(1)) := \Si^{\oplus} \, \Mc^{0}_{2k}
(\G(1)) \, ,
$$
the algebra of modular (resp. cusp) forms.
\smallskip

\noindent A richer and more interesting picture emerges when
lattices are replaced by pairs $ \, ( \Lba , \phi) \,$, with $
\Lba $ a lattice in $\mathbb{C}$, and
 $\displaystyle \phi : \frac{\mathbb{Q} \Lba}{\Lba}
\build\longrightarrow_{}^{\sim} \frac{\mathbb{Q}^2}{\mathbb{Z}^2}$
an isomorphism which is `unimodular' in a suitably defined sense,
cf. \cite[\S 17]{Mu}. The set $\, \Lc_\mathbb{A} \, $ of all such
pairs forms again a line bundle, only this time over a projective
limit of Riemann surfaces,
$$
 \mathbb{C}^{\ts} \, \ra \, \Lc_\mathbb{A} \, \ra \,
 {\FH}_\mathbb{A} := \quad \build\limproj_{N}^{} \G(N) \bash {\FH} \, ,
$$
where
$$ \G(N) \, = \, \left\{ \, \begin{pmatrix} a &b \\
 c &d  \end{pmatrix} \in \SL (2 , \mathbb{Z}) \, ;
\quad \begin{pmatrix} a &b \\
 c &d  \end{pmatrix} \, \equiv \, \begin{pmatrix} 1 &0 \\
 0 &1  \end{pmatrix} \quad \hbox{mod} \, N \, \right\} \,.
$$
and the projective system is ordered by divisibility. The
subscript $\mathbb{A}$ is justified by the adelic interpretation
of the fibration. Indeed, if one drops the `unimodularity'
condition on $\phi$, then the resulting space of pairs can be
canonically identified to the homogeneous space
\begin{equation} \label{adelgood}
 \GL (2,\mathbb{Q}) \bash \GL (2, \mathbb{A}) .
\end{equation}
The space $ \Lc_\mathbb{A}$ is the closure of $\GL (2,\mathbb{R})$
in this homogeneous space, and its structure of principal  bundle
comes from the natural inclusion of the multiplicative group
$\mathbb{C}^*$ in $\GL (2,\mathbb{R})$. One can also describe
${\FH}_\mathbb{A}$  cf. {\it loc. cit.} as a double quotient of
the form,
\begin{equation} \label{adel}
{\FH}_\mathbb{A}  \, \simeq \,
 \GL (2,\mathbb{Q}) \bash \GL (2, \mathbb{A})^0 / K_{\ify} \, Z_{\ify} \, ,
\end{equation}
\smallskip
where $ K_{\ify} \simeq S^{1}$, $Z_{\ify} \simeq \GL (1,
\mathbb{R})$ viewed as diagonal matrices in $\GL (2, \mathbb{R})$,
and the subgroup $\GL (2, \mathbb{A})^0 $ is defined by requiring
that the determinant belongs to $\mathbb{Q}^* \times \mathbb{R}^*
\subset \GL (1, \mathbb{A})$.
\medskip

\noindent Modular (resp. cusp) forms $\, f \,$ of weight $k$ with
respect to any  congruence subgroup  $ \G $ are defined in a
similar way as for $ \G(1) $, by requiring the holomorphicity
(resp. cuspidality) of $\, f|_{k} \, \g \,$ at ${\, \ify} \,$,
with respect to any local parameter $\, q^{1/m} \, $, for every
$\, \g \in \G(1) \,$. One obtains, for each $N \geq 1$, graded
algebras of {\it forms of level $\, N$}
$$
\Mc (\G(N)) := \Si^{\oplus}_{k \geq 1}
 \, \Mc_{k} (\G(N)) \, , \quad \hbox{resp.}
\quad \Mc^{0} (\G(N)) := \Si^{\oplus}_{k \geq 1} \, \Mc^{0}_{k}
(\G(N)) \, ,
$$
which can be assembled together into algebras of modular (resp.
cusp) forms of all levels:
$$
\Mc :=  \limind_{N \ra \ify} \, \Mc (\G(N)) \, , \quad
\hbox{resp.} \quad \Mc^{0} := \limind_{N \ra \ify} \, \Mc^0
(\G(N)) \, .
$$
\smallskip

\noindent The group
$$
G^+ (\mathbb{Q}) := \GL^{+} (2 , \mathbb{Q})  \, ,
$$
acts `sideways' on the tower defining the projective limit $\
{\FH}_\mathbb{A}$. Equivalently, one can view $G^+ (\mathbb{Q})$
as diagonally embedded in $  \GL (2 , \mathbb{A}_{f})^0 $, with
the latter acting on $\ {\FH}_\mathbb{A}$ by right translations
(cf. (\ref{adel})). Using the fact that one can factor any $\, \a
\in G^+ (\mathbb Q) \,$ in the form
\begin{equation} \label{bor}
\a = \g \cdot \b \, ,  \qquad \g \in \G (1) \, , \quad  \b \in B^+
(\mathbb Q) =
\left\{ \begin{pmatrix} a &b \\
 0 &d  \end{pmatrix} \in G^+ (\mathbb Q) \right\} \, ,
\end{equation}
it is easy to see that the action of $G^+ (\mathbb{Q})$ preserves
both holomorphicity and cuspidality. Thus, to this action of $G^+
(\mathbb{Q})$ corresponds a first algebra of `noncommutative
coordinates' on the ` transverse space' $ \, \Lc_\mathbb{A} / G^+
(\mathbb{Q}) \,$, namely the crossed-product algebra
$$ \Ac \, \equiv \,
\Ac_{G^+ (\mathbb{Q})} \, := \, \Mc \rtimes G^+ (\mathbb{Q}) \, ,
$$
containing as an ideal the crossed product
$$ \Ac^0 \, := \, \Mc^0 \rtimes G^+ (\mathbb{Q}) \, .
$$
The elements of $\Ac$  are finite sums of symbols of the form
$$ \sum \, f \, U_{\g}^* \, , \qquad \text{with} \quad f \in \Mc \, ,
\quad \g \in  G^+ (\mathbb{Q})\, ,
$$
and with the product given by the rule
\begin{eqnarray}
f \, U_{\a}^* \,\cdot  g \, U_{\b}^*  =\,( f\,\cdot g| \a)
 \;  U_{\b \,  \a}^* \, ,
\end{eqnarray}
where the `slash operation' is relative to the weight of $g$ as in
(\ref{act1}).
\medskip

\noindent On the other hand, to the right action of $  \GL (2 ,
\mathbb{A}_{f})^0 $ corresponds the inductive limit of the
following \textit{level N algebras}, $\Ac(N):= \Ac(\G(N))$,
obtained as crossed products of modular forms by Hecke operators
at level $N$, which we proceed now to describe.

\noindent More generally, for any congruence subgroup  $\G$ we
shall construct an algebra $\Ac(\G)$, which contains as
subalgebras both  the algebra of $\G$-modular forms $\Mc (\G)$ as
well as the Hecke ring at level $\G$ (comp.~\cite[Chap. 3]{Sh}),
without being in general generated by these subalgebras.
\bigskip

\noindent \textbf{Definition}.  Let $\G$ be  a congruence
subgroup. By a \textit{Hecke operator form of level} $\G$  we mean
a map
$$
 F: \G\bash G^+ (\mathbb{Q}) \ra \Mc \, ,
\qquad \G \a \mapsto F_{\a} \in \Mc  \, ,
$$
with \textit{finite support} and satisfying the \textit{covariance
condition}
\begin{equation} \label{tcov}
F_{\a \g} \, = \, F_{\a} \vert \g \, , \quad \fl \a \in G^+
(\mathbb{Q}) , \g \in  \G \, .
\end{equation}
A Hecke operator form $F$ of level $\G$ is called
\textit{cuspidal} if
\begin{equation} \label{cusp}
F_{\a} \in \Mc^0 \, , \quad \fl \a \in G^+ (\mathbb{Q}) \, .
\end{equation}
\medskip

\noindent Thus, by definition, $F$ is determined by the (nonzero)
values taken on a finite number of inequivalent representatives
$\{ \a_1 , \ldots , \a_r \}$ for double cosets in $\G\bash G^+
(\mathbb{Q})/\G \,$. It should be noted though that these values
$F_{\a_i } \in \Mc \,$ are not required to be modular forms of
level $\G$. Maps from $\G\bash G^+ (\mathbb{Q})/\G$ to $\Mc (\G)$
trivially fulfill equation (\ref{tcov}), but they do not exhaust
all its solutions. In fact, one has:
\bigskip

\begin{lemma} Let $f \in \Mc$ be any modular form, $\G$ a
congruence subgroup. There exists a Hecke operator form  $F$ of
level $\G$ such that $\, F_{\a} \, = \, f \, $ if and only if
\begin{equation*}
f \vert \g = f \, , \quad \fl \, \g \in \G \, \cap \a^{-1} \G \,
\a \, .
\end{equation*}
\end{lemma}
\medskip

\begin{proof} Assume $ \, f = F_{\a} \, $. Then for any
$\, \g \in \G \cap \a^{-1} \, \G \, \a \,$ one has $\, \a \g = \g'
\a \,$, for some $\,\g' \in \G \,$, and therefore $\, F_{\a \g} =
F_{\g' \a} = F_{\a} \,$. On the other hand, by (\ref{tcov}), $\,
F_{\a \g} = F_{\a} \vert \g \,$, which gives the required
invariance property of $f$.
\smallskip

\noindent Conversely, define $F$ to be $0$ on the complement of
the coset $ \, \G \, \a \, \G \,$, while
$$ F_{\g_1 \a \g_2} \, = \, f \vert \g_2 \, ,
\qquad \fl \, \g_1 \, , \g_2 \, \in \G \, .
$$
The definition is unambiguous, since if $ \, \g'_1 \a \g'_2 = \g_1
\a \g_2 \,$ then
$$
 {\g'}_1^{-1} \g_1 \a =
\a \g'_2 \g_2^{-1} \, , \quad \text{with}
 \quad \g'_2 \g_2^{-1} \in
 \G \, \cap \a^{-1} \G \, \a \, ,
$$
so that $\, f \vert \g'_2 = f \vert \g'_2 \g_2^{-1} \vert \g_2 = f
\vert \g_2 \,$. The condition (\ref{tcov}) is automatically
fulfilled.
\end{proof}
\bigskip

\noindent As a simple illustration let us consider the case $\G
=\G(1)$. The  double cosets in $\G(1)\bash G^+ (\mathbb{Q})/\G(1)
\,$ are
 represented by the elements $\a = r \,  \b_n$ where
 $r \in \GL^{+} (1 , \mathbb{Q})$ is the diagonal matrix
$  \begin{pmatrix} r &0 \\ 0 &r  \end{pmatrix} \,$, $n \in
\mathbb{N}$ and $ \, \b_n = \begin{pmatrix} n &0 \\ 0 &1
\end{pmatrix} \,$. Then
$$ \G(1) \cap \a^{-1} \, \G(1) \, \a \, = \, \G_0 (n) \,
:= \, \left\{ \, \begin{pmatrix} a &b \\
 c &d  \end{pmatrix} \in \G(1) \, ; \quad
 c \equiv \, 0  \quad (\text{mod} \, n) \, \right\} \,
$$
and any $\G_0 (n)$-modular form gives rise by the above recipe to
a Hecke operator form of level $\G(1)$. One can thus  view a
Hecke operator form of level $\G(1)$ as a  sequence with finite
support  of modular forms  $f_{n,r} \in \Mc(\G_0 (n))$. Note that
the action (\ref{act1}) on modular forms is trivial on the
component
 $ \GL^{+} (1 , \mathbb{R})$ of the center of $\GL(2, \mathbb{R})$
and therefore the label $r$ plays an unimportant role in the above
labelling. We keep it however, in order to retain the ability to
twist by characters of the center, for instance using the 
one-parameter group $\sigma_z, \, z\in \mathbb{C}$
 of automorphisms of  $\Ac(\G)$
given by
$$
(\sigma_z(F))_{\alpha}=\,\text{det}(\alpha)^z \, F_{\alpha}
$$
These are automorphisms for the algebra structure defined below in
Proposition 2.
 The weight of modular forms defines a natural grading
on $\Ac(\G)$, and one can see in the above example that the
algebra  $\Ac(\G(1))$ of Hecke operator forms of level $\G(1)$ is
non-trivial in weight $2$, while its subalgebra generated by
modular forms of level $1$ and standard Hecke operators has no
element of weight $2$.
\bigskip

\noindent Let us show  that the space $\Ac(\G)$ of Hecke
operator forms of level $\G$ is indeed an algebra.  
We shall refer to it as
the \textit{modular Hecke algebra of level} $\G$. This algebra
acts canonically on the space $\Mc (\G)$ of modular forms of level
$\G$, preserving the subspace
 $\, \Mc^0 (\G)$ of cusp forms.
\bigskip

\begin{proposition} Fix a congruence subgroup $\G$.
\begin{itemize}
\item[$1^0$.] With $\, F^1 , F^2 \in \Ac(\G) $, the operation
\begin{equation} \label{op}
(F^1 *  F^2)_{\a} \,:= \, \sum_{ \b \in \G \bash G^+ (\mathbb{Q})}
\, F^1_{\b} \cdot F^2_{\a \b^{-1}} \vert \b
\end{equation}
turns the vector space $\Ac(\G)$ of all Hecke operator forms of
level $\G$ into an associative algebra.
\item[$2^0$.] The equality
\begin{equation} \label{act}
F * f \, :=  \, \sum_{ \a \in \G \bash G^+ (\mathbb{Q})} \, F_{\a}
\cdot f \vert \a
\end{equation}
endows $\Mc  (\G)$ with the structure of an $\Ac(\G)$-module.
\item[$3^0$.] The cuspidal  Hecke operator forms of level
$\G$ form an ideal $\Ac^0(\G)$ of $\Ac(\G)$; the subspace
 $\, \Mc^0 (\G)$ is
a submodule of the $\Ac(\G)$-module $\Mc (\G)$, and $\Ac^0(\G)$
maps $\Mc (\G)$ to $\Mc^0 (\G)$.
\end{itemize}
\end{proposition}
\medskip

\begin{proof}  $1^0$. The right hand side of (\ref{op}) is
    well-defined, since the expression
$\,  F^1_{\b} \cdot F^2_{\a \b^{-1}} \vert \b \,$ is invariant
under the left multiplication of $\b$ by any $\g \in \G$.
Evidently, it defines a map with finite support from $\G\bash G^+
(\mathbb{Q}) $ to $ \Mc $. Also, for any $\, \g \in \G$,
\begin{eqnarray}
(F^1 * F^2)_{\a \g} \, &=& \,\sum_{ \b \in \, \G\bash G^+
(\mathbb{Q})} \,  F^1_{\b} \cdot F^2_{\a \g \b^{-1}} \vert \b \cr
\cr &=& \, \sum_{ \d \in \, \G\bash G^+ (\mathbb{Q})}
 F^1_{\d} \vert \g \cdot F^2_{\a \d^{-1}} \vert \d \vert \g \cr \cr
&=&  (F^1 * F^2)_{\a} \vert \g \, , \nonumber
\end{eqnarray}
so that $\, F^1 * F^2 \in \Ac(\G)$.

\noindent It is convenient to rewrite the product (\ref{op}) in
the form
\begin{equation} \label{op'}
(F^1 * F^2)_{\a} \, = \, {\sum}_{ \a_2 \, \a_1 = \a} \, F^1_{\a_1}
\cdot F^2_{\a_2} \vert \a_1 \, ,
\end{equation}
where the sum is taken over pairs $ (\a_1, \a_2 ) \in G^+
(\mathbb{Q}) $ such that $ \a_2 \, \a_1 = \a $, modulo the
equivalence relation
\begin{equation} \label{equiv}
(\a_1, \a_2 ) \, \sim \, (\g \a_1, \a_2 \g^{-1}) \, , \quad \g \in
\G \, .
\end{equation}
With this understood, the associativity of the product follows
from the fact that, irrespective of the placement of putative
parantheses,
$$
(F^1 * F^2 * F^3) _{\a} \, = \,  {\sum}_{ \a_3 \a_2 \a_1 = \a} \,
F^1_{\a_1} \cdot F^2_{\a_2} \vert \a_1
 \cdot F^3_{\a_3} \vert \a_2 \a_1 \, ,
$$
where the triples $ (\a_1, \a_2 , \a_3 ) \in G^+ (\mathbb{Q}) $
with $ \a_3 \a_2 \a_1 = \a $ are identified as above via the
equivalence
$$ (\a_1, \a_2 , \a_3 ) \, \sim \,
(\g \a_1, \d \a_2 \g^{-1} , \a_3 \d^{-1}) \, , \quad \g, \d \in \G
\, .
$$
\smallskip

\noindent $2^0$. To check that $ F * f \,$ is $\G$-modular, let $
\g \in \G \,$; then
$$
(F * f) \vert \g  =  \sum_{ \a \in \, \G\bash G^+ (\mathbb{Q})} \,
F_{\a} \vert \g \cdot f \vert \a \g =
 \sum_{ \a \in \, \G\bash G^+ (\mathbb{Q})}
 \, F_{\a \g} \cdot f \vert \a \g = F * f \, .
$$
The above proof of associativity also applies to show that we do
have an action of $\Ac(\G)$ on the vector space $\Mc (\G)$.
\smallskip

\noindent Finally, $3^0 \,$ follows from the very definitions,
using the ``Iwasawa decomposition'' (\ref{bor}).
\end{proof}
\bigskip

 \noindent Let $\Hc(\G)$ be the standard Hecke algebra of functions
with finite support on double cosets in $\G\bash G^+
(\mathbb{Q})/\G \,$. It is canonically isomorphic to the
subalgebra of weight $0$ elements of  $\Ac(\G)$ by the isomorphism
$$
j: \Hc(\G) \rightarrow \Ac(\G)\, , \,\, j(h)_{\alpha}:= h(\alpha)
 \quad h \in  \Hc(\G)
$$
The augmentation $\epsilon$ of the graded algebra  $\Ac(\G)$ gives
a homomorphism $\epsilon:\Ac(\G) \rightarrow \Hc(\G)$ such that
$\epsilon \circ j \,= \,\text{Id}$.
\bigskip

\noindent For $\G=\G(N)$ a Hecke operator form of level $\G$ is
simply called ``of level $N$''. It  gives rise for each multiple
$N'$ of $N$ to a map
$$
 F': \G(N')\bash G^+ (\mathbb{Q}) \ra \Mc \, ,
\qquad F' := F \, \circ \, p ,
$$
where $p$ is the projection,
$$
\G(N')\bash G^+ (\mathbb{Q}) \rightarrow \G(N)\bash G^+
(\mathbb{Q})
$$
The map $F'$ has finite support and the covariance condition
\begin{equation} \label{tcov1}
F'_{\a \g} \, = \, F'_{\a} \vert \g \, , \quad \fl \a \in G^+
(\mathbb{Q}) , \g \in \G(N') \, .
\end{equation}
is still fulfilled. Let $r(N,N')$ be the order of $\G(N')\bash
\G(N)$, we then define a map $\rho_{N,N'}$ from $\Ac(N)$ to
$\Ac(N')$ by
$$
\rho(F) := \frac{1}{r(N,N')} F \, \circ \, p
$$
for every Hecke operator form $F$ of level $N \,$. This defines a
homomorphism from $\Ac(N)$ to $\Ac(N')$ and allows to take the
inductive limit,
$$
\Ac(\infty) :=  \limind_{N \ra \ify} \, \Ac (N) \, ,
$$
which in turn can be interpreted as the crossed product of $\Mc$
by $  \GL (2 , \mathbb{A}_{f})^0 $.

\medskip
\section{The Hopf algebra $\Hc_1$ and its Hopf cyclic classes}

In this section we describe the Hopf algebra $\Hc_1$ and its Hopf
cyclic cocycles corresponding respectively to the Godbillon-Vey
class, to the Schwarzian derivative, and to the transverse
fundamental class.
\medskip

\noindent We begin by recalling the presentation of the Hopf
algebra $\Hc_1$ (cf. \cite{CM1}). As an algebra, it coincides with
the universal enveloping algebra of the Lie algebra with basis $\{
X,Y,\d_n \, ; n \geq 1 \}$ and brackets
\begin{equation} \label{pres}
[Y,X] = X \, , \, [Y , \d_n ] = n \, \d_n \, , \, [X,\d_n] =
\d_{n+1} \, , \, [\d_k , \d_{\ell}] = 0 \, , \quad n , k , \ell
\geq 1 \, .
\end{equation}
As a Hopf algebra, the coproduct $\, \D : \Hc_1 \ra \Hc_1 \ot
\Hc_1 \,$ is determined by
\begin{eqnarray}
\D \,  Y = Y \ot 1 + 1 \ot Y \, , \quad
\D \,  X &=& X \ot 1 + 1 \ot X + \d_1 \ot Y \nonumber \\
\D \,  \d_1 &=& \d_1 \ot 1 + 1 \ot \d_1
\end{eqnarray}
and the multiplicativity property
\begin{equation}
\D (h^1 \, h^2) = \D h^1 \cdot \D h^2 \, , \quad h^1 , h^2 \in
\Hc_1 \, \, ;
\end{equation}
the antipode is determined by
\begin{equation}
S(Y) = -Y \, , \, S(X) = -X + \d_1 Y \, , \, S(\d_1) = - \d_1
\end{equation}
and the anti-isomorphism property
\begin{equation}
S (h^1 \, h^2) = S(h^2) \, S (h^1) \, , \quad h^1 , h^2 \in \Hc_1
\, \, ;
\end{equation}
finally, the counit is
\begin{equation}
\ve (h) = \hbox{constant term of} \quad h \in \Hc_1 \, .
\end{equation}

\noindent The role of $\Hc_1$ as symmetry in tranverse geometry
comes from its natural action on crossed products (\cite{CM1}).
Given a one-dimensional manifold $M^1$ and a discrete subgroup $\G
\sbs \Diff^+ (M^1)$, $\, \Hc_1$ acts on the crossed product
algebra
$$\, \Ac_{\G} = C_c^{\ify} (J_+^1 (M^1)) \rtimes \G \, ,
$$
by a Hopf action, where $J_+^1 (M^1)$ is the oriented $1$-jet
bundle over $M^1$. We use the coordinates in $J_+^1 (M^1)$ given
by the Taylor expansion,
$$
j(s) = y + s \, y_1 +  \cdots \, , \qquad y_1 > 0 \, ,
$$
and let diffeomorphisms act in the obvious functorial manner on
the $1$-jets,
$$
\vp (y, y_1) = (\vp (y) , \, \vp' (y) \cdot y_1) \, .
$$
The action of $\Hc_1$ is then given
 as follows:
\begin{eqnarray}
Y(f U_{\vp}^*) = y_1 \, \frac{\partial f}{\partial y_1} \,
U_{\vp}^* \, &,& \qquad X(f U_{\vp}^*) = y_1 \, \frac{\partial
f}{\partial y} \, U_{\vp}^*
\, , \\
\, \d_n (f U_{\vp}^*) &=& y_1^n \,  \frac{d^{n}}{d y^{n}}
\left(\log \frac{d\vp}{dy}\right) \, f U_{\vp}^* \, \, ,
\label{dn}
\end{eqnarray}
where we have identified $\, J_+^1 (M^1) \simeq M^1 \times
\mathbb{R}^+ \, $ and denoted by  $\, (y, y_1) \, $ the
coordinates on the latter.
\medskip

\noindent The volume form $\displaystyle \frac{dy \wdg
dy_1}{y_1^2} \, $ on $ \, J_+^1 (M^1) \,$ is invariant under $\,
\Diff^+ (M^1) \,$ and gives rise to the following trace $\, \tau :
\Ac_{\G} \ra \mathbb{C} $,
\begin{equation}
\tau (f U_{\vp}^*) \, = \, \begin{cases}
 \int_{J_+^1 (M^1)} f(y,y_1) \, \frac{dy \wdg dy_1}{y_1^2}
 &\text{ if }  \vp = 1 \, ,  \\
0 &\text{ if }  \vp \ne 1 \, .
\end{cases}
\end{equation}
This trace is $\nu$-invariant with respect to the action $\, \Hc_1
\ot \Ac_{\G} \ra \Ac_{\G} \,$ and with the modular character $\,
\nu \in \Hc_1^* \,$, determined by
$$\nu (Y) = 1, \quad \nu (X) = 0 , \quad \nu (\d_n) = 0 \, ;
$$
the invariance property is given by the identity
\begin{equation}
\tau (h(a)) = \nu (h) \, \tau (a) \, , \qquad \fl \quad h \in
\Hc_{1} \, .
\end{equation}
\smallskip

\noindent The fact that
$$\, S^2 \ne \Id \, ,
$$
is automatically corrected by twisting with $\nu$. Indeed, $\wt S
= \nu * S$ satisfies the involutive property
\begin{equation} \label{inv1}
\wt S^2 = \Id \, .
\end{equation}
One has
\begin{equation} \label{tildes}
\wt S(\delta_1)\, = \, -\delta_1 \,, \quad \wt S(Y) = - Y +1 \, ,
\quad \wt S(X)= -X + \delta_1 Y \, .
\end{equation}
Equation (\ref{inv1}) shows that the pair ($\nu,1$) given by the
character $\nu$ of $\Hc_1$ and the group-like element $1 \in
\Hc_1$ is a modular pair in involution, which thus allows us to
define the cyclic cohomology $ HC_{\Hopf}^* \,(\Hc_1)$ (cf.
\cite{CM1}, \cite{CM2}). We refer to Appendix A for the detailed
definition of the Hopf cyclic cohomology.
\medskip

\noindent   The assignment
\begin{equation}
\chi_{\tau} (h^{1} \ot \ldots \ot h^{n}) (a^0 , \ldots , a^{n}) \,
= \, \tau (a^0 \, h^{1} (a^{1}) \ldots h^{n} (a^{n})) \, ,
\end{equation}
where $ \, h^{1}, \ldots , h^{n} \, \in \Hc_{1} \,$ and $\, a^{0},
a^{1}, \ldots , a^{n} \, \in \Ac_{\G} \,$, induces a
characteristic homomorphism
$$
\chi_{\tau}^* \, : \,  HC_{\Hopf}^* \,(\Hc_1) \, \ra \, HC^* \,
(\Ac_{\G}) \, .
$$
\medskip

\noindent In \cite{CM1} we have constructed an isomorphism
$$ \kappa_{1}^{*} \, : \, H^{*} ({\Fa}_1 , \mathbb{C})
\build\longrightarrow_{}^{\simeq}P HC_{\Hopf}^{*} \, (\Hc_1) \,,
$$
between the Gelfand-Fuchs cohomology of the Lie algebra of formal
vector fields on $\mathbb{R}^{1}$ and the periodic cyclic
cohomology of the Hopf algebra $\Hc_1$.
\bigskip

\begin{proposition} \label{gv2}
The element $\, \d_{1} \in \Hc_{1} \,$ is a Hopf cyclic cocycle,
which gives a nontrivial class
$$ [\d_{1}] \, \in \, HC_{\Hopf}^{1} \, (\Hc_1) \, .
$$
Moreover, $ \, [\d_{1}] \,$ is a generator for $ \,
PHC_{\Hopf}^{\rm odd} \, (\Hc_1) \,$ and corresponds to the
Godbillon-Vey class in the isomorphism $\kappa_{1}^{*}$  with the
Gelfand-Fuchs cohomology.
\end{proposition}

\begin{proof}
 Indeed, the fact that $\d_{1}$ is a $1$-cocycle is easy to check:
$$
b (\d_1) \, = \, 1 \ot \d_1 - \D \d_1 + \d_1 \ot 1\,  = \, 0 \, ,
$$
while
$$
\tau_1 (\d_1) \, = \, \wt S (\d_1) \, = \, S (\d_1) \, = \, - \d_1
\, .
$$
On the other hand, its image under the above characteristic map,
$$
\chi_{\tau}^* \,([\d_{1}]) \, \in \, HC^{1} \, (\Ac_{\G}) \, ,
$$
is precisely the anabelian $1$-trace of \cite{C1} (cf. also
\cite[III. 6. $\g$]{C2}), and the latter is known to give a
nontrivial class on the transverse frame bundle to codimension $1$
foliations. The remaining statement is proved in Appendix B
(Proposition \ref{gv1}).
\end{proof}
\bigskip

\noindent We shall now describe another Hopf cyclic $1$-cocycle
which corresponds to the Schwarzian derivative
 $\{y \, ;x\} \,$, whose expression is
\begin{equation}
\{y \, ;x\} \, := \, \frac{d^2}{dx^2} \left( \log \frac{dy}{dx}
\right) \, - \, \frac{1}{2} \left( \frac{d}{dx} \left( \log
\frac{dy}{dx} \right) \right)^2 \, .
\end{equation}
\medskip

\begin{proposition} \label{sc1}
The element $\, \d'_{2} := \d_2- \frac{1}{2}\d_1^2 \, \in \Hc_{1}
\,$ is a Hopf cyclic cocycle, whose  action on the crossed product
algebra $ \Ac_{\G} = C_c^{\ify} (J_+^1 (M^1)) \rtimes \G \, $ is
given by the Schwarzian derivative
$$
\d'_2 (f U_{\vp}^*) = y_1^2 \, \{\vp(y) \, ;y\} \, f U_{\vp}^*
$$
and whose class
$$  [\d'_{2}] \, \in \, HC_{\Hopf}^{1} \, (\Hc_1) \,
$$
is equal to $ B(c)$, where $c$ is the following Hochschild
$2$-cocycle,
$$ c \, := \, \delta_1 \ot X + \frac{1}{2} \, \delta_1^2 \ot Y \, .
$$
\end{proposition}
\begin{proof} We shall give the detailed computation in order
 to illustrate the $(b,B)$ bicomplex for Hopf cyclic cohomology.
Let us compute $b (c)$. One has
\begin{eqnarray}
b(\delta_1 \ot X) = 1 \ot \delta_1 \ot X - (\delta_1 \ot 1 + 1 \ot
\delta_1) \ot X + \\ \nonumber \delta_1 \ot(X \ot 1 + 1 \ot X +
\delta_1 \ot Y) - \delta_1 \ot X \ot 1  = \delta_1 \ot \delta_1
\ot Y
\end{eqnarray}
Also
\begin{eqnarray}
 b(\delta_1^2 \ot Y)= 1 \ot \delta_1^2 \ot Y - (\delta_1^2 \ot 1 +
2 \delta_1 \ot \delta_1 + 1 \ot \delta_1^2) \ot Y + \\ \nonumber
\delta_1^2 \ot ( Y \ot 1 + 1 \ot Y) - \delta_1^2 \ot Y \ot 1 = -2
\delta_1 \ot \delta_1 \ot Y
\end{eqnarray}
This shows that
$$
    b(c) \, = \, 0 \, ,
$$
so that $c$ is a Hochschild cocycle.
\medskip

\noindent Let us now compute $B(c)$. First, we recall that
$$
B_0 (h^1 \ot h^2) = \wt S (h^1) h^2 \, .
$$
Since $\wt S(\delta_1)\, = \, -\delta_1 \,$, one then has
$$
B_0 (c) = - \delta_1 \, X + \frac{1}{2}  \delta_1^2 \, Y \, .
$$
Since $\, \wt S(Y) = - Y +1 \,$ and $\, \wt S(X)= -X + \delta_1 Y
\,$, it follows that
\begin{eqnarray}
\wt S(B_0c)  &=&  \wt S(X) \delta_1 + \frac{1}{2}  \wt S(Y)
\delta_1^2  \\ \nonumber &=&(- X + \delta_1 Y)\delta_1 +
\frac{1}{2}  (- Y +1) \delta_1^2  \\ \nonumber
 &=&-X \delta_1 + \delta_1^2 Y + \delta_1^2 - \frac{1}{2}  (\delta_1^2 Y +
\delta_1^2)  \\ \nonumber
  &=&- X \delta_1 + \frac{1}{2}  \delta_1^2 Y + \frac{1}{2}  \delta_1^2 \, .
\end{eqnarray}
Therefore,
\begin{eqnarray}
B(c)  \, = \, B_0c - \wt S(B_0c)  \, &=& \,
\\ \nonumber
- \delta_1 \, X + \frac{1}{2}  \delta_1^2 \, Y - (-X \delta_1 +
\frac{1}{2}  \delta_1^2 Y + \frac{1}{2}  \delta_1^2) \,
 \, &=& \, \delta'_2 \, .
\end{eqnarray}
which shows that the class of $\d'_2$ is trivial in the periodic
cyclic cohomology $PHC_{\Hopf}^{*}(\Hc_1)$.
\end{proof}
\medskip

\noindent   We conclude this section by noting that while the
class of the unit constant $\, [1] \in HC_{\Hopf}^{0}(\Hc_1) \,$
is trivial in the periodic theory  (since $ \, B (Y) = 1 \,$), a
generator of  $PHC_{\Hopf}^{\rm even}(\Hc_1)$ is the class of the
cyclic $2$-cocycle
\begin{equation} \label{fund}
F := X \ot Y - Y \ot X - \d_1 \,  Y \ot Y \, ,
\end{equation}
which in the foliation context represents the `transverse
fundamental class'.

\medskip

\section{Hopf Symmetry of Modular Hecke  Algebras}

We shall define and analyze  in this section the natural action of
the Hopf algebra $\, \Hc_1 \,$ on the algebras $\, \Ac_{G^+
(\mathbb{Q})}  \,$ and $\Ac(\G)$, for any congruence subgroup
$\G$. We shall then formulate in precise terms the geometric
analogy with the action of $\, \Hc_1 \,$ on the crossed products
of the polynomial functions on the frame bundle of $S^1$ by
discrete subgroups of $\Diff(S^1)$.
\medskip

\noindent In order to define the action of the generator $X$, we
shall use the most natural derivation of the algebra of modular
forms, which is in fact a classical one, originating in the work
of Ramanujan (cf. \cite{Ra}, \cite{Se}). It is the operator of
degree 2 on the algebra of modular forms given by
\begin{equation} \label{ram}
X := \frac{1}{2  \pi i} \frac{d}{dz} - \frac{1}{12  \pi i} \,
\frac{d}{dz} (\log \D) \cdot Y = \frac{1}{2  \pi i} \,
\frac{d}{dz} - \frac{1}{2  \pi i} \, \frac{d}{dz} (\log \eta^4)
\cdot Y \, ,
\end{equation}
where
\begin{equation}
\D (z) \, = \, (2  \pi)^{12} \, \eta^{24} (z) \, = \, (2
\pi)^{12} \, q \, \prod_{n=1}^{\infty} (1 - q^n)^{24} \, , \quad q
= e^{2\pi i z}
\end{equation}
is the discriminant cusp form of weight 12, expressed in terms of
the Dedekind $ \eta$-function, while $\, Y \,$ stands for the
grading operator
\begin{equation}
Y(f) = \frac{k}{2} \cdot f \, , \quad \fl \, f \in \Mc_{k} \, .
\end{equation}
The derivation $\, X \,$ provides the analogue of a connection,
given by means of its horizontal vector field, while the operator
$Y$ serves as the analogue of the fundamental vector field.
Indeed, one has
$$ [Y, X] \, = \, X \, .
$$
\noindent As a matter of fact, if we interpret the modular forms
of weight $k$ as sections of the $\frac{k}{2}$ power of the line
bundle of $1$-forms, via the identification
\begin{equation} \label{highdiff}
f \mapsto f \,\cdot \,(2 \pi i\, dz)^{\frac{k}{2}}\, , \quad \fl
\, f \in \Mc_{k} \, ,
\end{equation}
$X$ effectively defines a connection $\nabla$, uniquely determined
by the property
\begin{equation} \label{vanish}
\nabla(\eta^{2k}dz^{\frac{k}{2}}) \, = \, 0 , \qquad k \in 2
{\mathbb N} \, \, ;
\end{equation}
explicitly,
\begin{equation} \label{connect}
\nabla(f \cdot (2\pi i\, dz)^{\frac{k}{2}}) \, = \, X(f) \cdot (2
\pi i\, dz)^{\frac{k}{2} + 1} \, , \quad \fl \, f \in \Mc_{k} \, .
\end{equation}
With this understood, we can now state the following:
\bigskip

\begin{lemma} \label{mu}
For any $\, \g \in G^+ (\mathbb{Q})\,$ and $\, f \in \Mc_{k} \,$,
one has
$$
    \left( X(f|_{k} \, {\g^{-1}}) \right)|_{k + 2} \, \g \,
= \, X(f) \, - \, \mu_{\g} \cdot Y(f) \, ,
$$
where
\begin{equation} \label{dmu}
\mu_{\g} \, (z) \, = \, \frac{1}{12  \pi i} \, \frac{d}{dz} \log
\frac{\D | \g}{\D} \, = \, \frac{1}{2  \pi i} \, \frac{d}{dz} \log
\frac{\eta^4 | \g}{\eta^4} \, .
\end{equation}
\end{lemma}
\smallskip

\begin{proof} This follows from (\ref{connect}), with $\mu_{\g}$
accounting for the lack of invariance of the section $\eta^{4}\,
dz$.

\end{proof}

\bigskip

\noindent {\bf Remark 1. \,} Alternatively, one can proceed by
direct calculation and use (\ref{ram}) to obtain,
for any $\, \g = \begin{pmatrix} a &b\\
 c &d  \end{pmatrix} \in G^+ (\mathbb{Q})\,$,
\begin{equation*}
 \mu_{\g} \, (z) \, = \, \frac{1}{12  \pi i} \left(
\det (\g) (cz + d)^{-2} \, \frac{\D'}{\D} \left( \frac{az+b}{cz+d}
\right) - \frac{\D'}{\D} \, (z) - \frac{12c}{cz+d} \right) \, .
\end{equation*}
The right hand side is the explicit form of that in (\ref{dmu}).
To interpret it, we recall that (cf. e.g. \cite[III, \S 2]{Sch})
\begin{equation} \label{rmu}
\frac{\D^{'}(z)}{\D(z)}\, = \, - \frac{6}{\pi i} G^*_2 (z) \, ,
\end{equation}
where
\begin{equation} \label{G*}
 G^*_2 (z) \, = \, 2 \zeta (2) \, + \,
2 \sum_{m \geq 1} \, \sum_{n \in \mathbb{Z}} \frac{1}{(mz + n)^2}
\, = \, \frac{\pi^2}{3} \, - \, 8 \pi^2 \sum_{m, n \geq 1} m e^{2
\pi i m n z}
\end{equation}
is the holomorphic Eisenstein series of weight 2. The latter fails
to be modular, and satisfies the transformation formula
\begin{equation} \label{g*}
  G^*_2 | \a \, ( z) \, = \,
G^*_2(z) \, - \,\frac{2 \pi i \, c}{cz+d}  \, , \qquad \fl \, \a =
\, \begin{pmatrix} a &b\\
 c &d  \end{pmatrix} \in \G(1) \, .
\end{equation}
Thus, $\, \mu_{\a}\equiv 0 \,$ when $\, \a \in \SL (2, \mathbb{Z})
\,$,
while for $ \g = \begin{pmatrix} a &b\\
 c &d  \end{pmatrix} \in G^+ (\mathbb{Q})\ $
\begin{equation} \label{nmu}
 \mu_{\g} \, (z) = \frac{1}{2 \, \pi^2}
 \left(  G^*_2 | \g \, ( z) -  G^*_2 ( z)
+ \frac{2  \pi i \, c}{cz+d} \right)
\end{equation}
measures the failure of the gauge transformation $\g$ of
preserving the `connection'.
\bigskip

\noindent By its very definition (\ref{dmu}), $\, \mu \,$ takes
values in the space of weight $2$ modular forms $\, \Mc_2 \,$. It
is not difficult to see that in fact the range of $\, \mu \,$ is
contained in the canonical complement $\, \Ec_2 \sbs \Mc_2 \,$ to
the  subspace of cuspidal modular forms $\, \Mc_2^0 \,$  ,
generated by Eisenstein series. Although we shall later prove a
more precise result, let us record for now the following:
\bigskip

\begin{lemma} \label{eis}
For any $\, \g \in G^+ (\mathbb{Q}) \,$,
 $\, \mu_{\g} \,$ is an
 Eisenstein series of weight 2.
\end{lemma}
\smallskip

\begin{proof} From (\ref{dmu})
it follows that for any $ \g_{1}\, , \g_{2} \in G^+ (\mathbb{Q})\
$
\begin{equation} \label{gmu}
\mu_{\g_{1}\cdot \g_{2}} \, = \, \mu_{\g_{1}}| \g_{2} \, + \,
\mu_{\g_{2}} \, \, ;
\end{equation}
in particular,
\begin{equation} \label{zmu}
\mu_{\g_{1}\cdot \g_{2}} \, = \, \mu_{\g_{2}} \, , \quad \fl \,
\g_1 \in \G(1) \, .
\end{equation}

\noindent Using the decomposition
$$ G^+ (\mathbb{Q}) \, = \, \G (1) \cdot T^+ (\mathbb{Q}) \cdot \G (1) \,
\, , \quad \hbox{where} \quad T^+ (\mathbb{Q}) =
\left\{ \begin{pmatrix} a &0\\
 0 &d  \end{pmatrix} \in G^+ (\mathbb{Q}) \right\} \, ,
$$
one can write any $\, \g \in G^+ (\mathbb{Q}) \,$ as a product
$$
\g  \, = \, \a_1 \cdot \d \cdot \a_2 \, ,  \qquad \a_1 , \a_2 \in
\G (1) \quad \hbox{and} \quad \d \in T^+ (\mathbb{Q})  \, \, ;
$$
when $\, \g \,$ has integer entries, so does $\, \d$. Applying
(\ref{gmu}) and the fact that $\, \mu \,$ vanishes on $\, \G (1)
\,$ one obtains:
\begin{equation} \label{gd}
\mu_{\g} \, = \, \mu_{\a_1} | \d \cdot \a_2 \, + \, \mu_{\d \cdot
\a_2 } \, = \, \mu_{\d} | \a_2 \, + \, \mu_{\a_2} \, = \, \mu_{\d}
| \a_2 \, .
\end{equation}
Thus, the range of $\, \mu \,$ is spanned by elements of elements
of the form
\begin{equation} \label{damu}
 \mu_{\d} | \a  \, = \, \frac{1}{2 \, \pi^2}  \cdot
 (  G^*_2 \, | \d \, -
 \, G^*_2 ) \, | \a  \, ,
 \quad \d  \in T^+ (\mathbb{Q}) \quad \hbox{and} \quad
\a \in \G(1) \, .
\end{equation}
The fact that the range of $\, \mu \,$ consists of weight $2$
Eisenstein series follows now from the identity, cf. \cite[VII, \S
3.5]{Sch},
\begin{equation} \label{em}
 N \, G^*_2 (N z) \, - \, G^*_2 (z) \, = \, N^{-1}
 \sum_{k=1}^{N-1} \, \wp_{(0,\frac{k}{N})} \, (z) \, ,
 \qquad N \geq 2 \,
\end{equation}
where
$     \wp_{\ab} (z) \,$
is the ${\ab}-$division value of the Weierstrass $\wp$-function
and the collection of functions
\begin{equation*}
 \left\{ \wp_{\ab} \, \, ; \,
\, {\ab} \in \left( \frac{1}{N} \mathbb Z \slash \mathbb Z
\right)^2 \setminus {\0b} \right\}
\end{equation*}
generates the space of weight $2$ Eisenstein series of level $N$.
\end{proof}
\bigskip

\noindent We now analyze the action of the Hopf algebra $\, \Hc_1
\,$ on the algebra $\, \Ac_{G^+ (\mathbb{Q})}  \,$, the
corresponding results for $\Ac(\G)$ for any congruence subgroup
$\G$ will be discussed afterwards.

\bigskip

\begin{proposition} \label{cor}
There is a unique Hopf action of $\Hc_1$ on $\, \Ac_{G^+
(\mathbb{Q})} \,$, determined by
\begin{eqnarray}
X (f \, U_{\g}^*) = X(f)\, U_{\g}^* \, , \qquad
Y (f \, U_{\g}^*) = Y(f)\, U_{\g}^* \, , \nonumber \\
\hbox{and} \qquad \d_1 (f \, U_{\g}^*) = \mu_{\g} \cdot f \,
U_{\g}^* \, . \qquad \qquad \qquad \nonumber
\end{eqnarray}
\end{proposition}
\smallskip

\begin{proof} Indeed, one has by construction
$$ [Y, X] \, = \, X \, , \quad [Y,  \d_{1}] \, = \,  \d_{1}
$$
The action of $\d_n$ is uniquely determined by induction using
$[X, \, \d_n]=\, \d_{n+1}$, which gives
\begin{equation} \label{deltan}
 \d_n (f \, U_{\g}^*) = X^{n-1}(\mu_{\g}) \cdot f \, U_{\g}^* \, .
\end{equation}
The relations $[\d_k, \, \d_n]=0$ follow easily. Since as an
algebra $\Hc_1$ is the envelopping algebra of the Lie algebra with
presentation
$$ [Y, X] \, = \, X  ,\, [Y,  \d_{n}] \, = \, n \d_{n}
, [X, \, \d_n]=\, \d_{n+1}, \,[\d_k, \, \d_n]=0,
$$
one gets the required action of $\Hc_1$. Let us check its
compatibility with the coproduct. For any $a^1 , a^2 \in \Ac_{G^+
(\mathbb{Q})}  \,$ one checks that,
$$
Y(a^1 a^2 ) \, = \, Y (a^1) \, a^2 + a^1 \, Y (a^2) \, , \qquad
\qquad \qquad
$$
$$
 X (a^1 a^2) \, = \, X (a^1) \, a^2 +
a^1 \, X (a^2) + \d_{1} (a^1) \, Y (a^2) \, ,
$$
$$
 \d_{1} (a^1 a^2) \, = \, \d_{1} (a^1) \, a^2 +
a^1 \, \d_{1} (a^2) \, ,
$$
which is enough to prove the required compatibility since $X, Y,
\d_{1}$ generate $\Hc_1$ as an algebra.
\end{proof}
\bigskip

\noindent We shall now show that, after a change of coordinates,
the action of $ \Hc_1$ on $\, \Ac_{G^+ (\mathbb{Q})} \,$ in the
new coordinates becomes strikingly similar to the action of $
\Hc_1$ on the crossed product algebra $\, \Ac_{\G} = C_c^{\ify}
(J_+^1 (M^1)) \rtimes \G \,$ (cf. \S 1). To this end, let us
introduce the primitive of the differential form $ \frac{2 \pi
i}{6}\eta^4 dz$, which is, in classical terminology, an integral
of the $1$st kind for $\G (6)$. Thus,
\begin{equation} \label{Z}
Z (z) \, := \,\frac{2 \pi i}{6} \int_{i \infty}^z \eta^4 dz \, ,
\end{equation}
and
\begin{equation} \label{newZ}
    dZ \, := \, \frac{1}{6} \, \eta^4 \, \frac{dq}{q} \,
    = \,\frac{2 \pi i}{6}
 \,  \eta^4 \, dz \, .
\end{equation}
Let $\chi$ be the unique character of $\G (1)$ such that
\begin{equation} \label{car}
    \chi (\begin{pmatrix} 1 &1\\
 0 &1  \end{pmatrix})\,=\, e^{\frac{2 \pi i}{6}}\,.
\end{equation}
One has by construction
\begin{equation*}
    dZ |_{0} {\g}
 \, = \, \chi (\g)\, dZ \quad \forall \g \in \G (1)
\end{equation*}
The kernel of $\chi$ is the commutator subgroup $\G'(1)$ which is
a free group on the two generators,
\begin{equation} \label{gen}
    \g_1=\begin{pmatrix} 2 &1\\
 1 &1  \end{pmatrix} \, ,\quad\,  \g_2=\begin{pmatrix} 1 &1\\
 1 &2  \end{pmatrix}
\end{equation}
For $\g \in \G'(1)$ one has $dZ |_{0} {\g}=dZ$ and $Z |_{0}
{\g}-\,Z = L(\g)$, where $L$ is the unique homomorphism from the
free group $\G'(1)$ to the additive group $\mathbb{C}$ such that,
\begin{equation} \label{lll}
 L(\g_1)=\, L_0 \, \, e^{\frac{2 \pi i}{3}}\, , \quad
L(\g_2)=\, L_0 \, \, e^{\frac{2 \pi i}{6}}
\end{equation}
where $2\,L_0 \sim 1.402182....$ is one third of the real half period
of the elliptic integral $\int\,\frac{dx}{\sqrt{x^3+1}}$. This is
easily seen by evaluating $Z |_{0} {\g_k}(i \infty)$.
 It follows in particular that $Z$ gives an isomorphism
of the elliptic curve $\G'(1) \bash H^* $ with the quotient
$\mathbb{C} / \Lba \,$ where $\Lba$ is the cubic lattice which is
the range of $L$. Extending $L$ to $\G (1)$ we obtain the
following,
\begin{equation} \label{chi}
    Z |_{0} {\g}
 \, = \, \chi (\g)\, Z \, + L(\g)
\quad \forall \g \in \G (1)
\end{equation}
We shall now show that the variable $Z \in \mathbb{C} / \Lba \,$
plays exactly the same role as the variable $\theta
 \in \mathbb{R} /2\,\pi\, \mathbb{Z} \,$ of the circle case.

\noindent Using formula  (\ref{connect}), one has
\begin{equation} \label{Xf}
X(f) \cdot (2  \pi i \, dz)^2 \, = \, d(\frac{ 2  \pi i \,f
\,dz}{dZ})\,dZ \, , \qquad \fl f \in \Mc_2 \, .
\end{equation}
Since
$$
2  \pi i \,\mu_{\g} \,dz= d \log \frac{\eta^4 | \g}{\eta^4} =
\frac{d}{dZ} \left( \log \frac{d(Z |_{0} {\g}) }{dZ} \right) dZ \,
,
$$
it follows from (\ref{Xf}) that
\begin{equation} \label{Xmu}
 X(\mu_{\g}) \cdot (2 \pi i \, dz)^2 \, = \, (\frac{d}{dZ})^2
 \left( \log \frac{d(Z |_{0} {\g})}{dZ} \right) (dZ)^2 \, .
\end{equation}

Iterating the above calculation one obtains for the higher
horizontal derivatives of $ \mu_{\g}$ the expresssion
\begin{equation}  \label{highmu}
 X^{n-1}(\mu_{\g}) \cdot (2  \pi i \, dz)^n \, = \,( \frac{d}{dZ})^n
 \left( \log \frac{d(Z |_{0} {\g})}{dZ} \right) (dZ)^n \, ,
 \quad \fl n \in {\mathbb N} \,
\end{equation}
\smallskip
so that we can state the following exact analogue of (\ref{dn}).
\medskip

\begin{proposition} \label{gend}
    The generalized cocycles  $\d_n \in \Hc_1$
act on $\, \Ac_{G^+ (\mathbb{Q})} \,$ by the formulae:
\begin{equation}
\d_n (f  \, U_{\g}^* ) \, = \, (\frac{d}{dZ})^n
 \left( \log \frac{d(Z |_{0} {\g}}{dZ} \right) (dZ)^n \,
  f  \, U_{\g}^* \, .
\end{equation}
\end{proposition}
\smallskip
\noindent In the above statement, we have used the identification
(\ref{highdiff}) between modular forms and higher differentials,
so that the modular form $f$ of weight $k$  stands here for the
higher differential $f \cdot (2  \pi i \, dz)^{\frac{k}{2}}$.
\bigskip

\noindent In the context of codimension $1$ foliations discussed
in Section 2, we showed (cf. Proposition \ref{sc1}) that the
$1$-cocycle
$$
\d_2' \, = \, \d_2 -\frac{1}{2} \d_1^2 \in \Hc_1 \, ,
$$
corresponds to the Schwarzian derivative. Here, using Proposition
\ref{cor}, we obtain:
\medskip

\begin{corollary} \label{sch}
The equality
\begin{equation}
\d_2' (f \, U_{\g}^*) = \left(X(\mu_{\g}) - \frac{1}{2}\,
\mu_{\g}^2 \right)
 \cdot f \, U_{\g}^*
\end{equation}
defines a derivation of the algebra $\, \Ac_{G^+ (\mathbb{Q})}
\,$.
\end{corollary}
\smallskip

\begin{proof} Indeed, one easily checks
 that the action of $\d_2'$ on $\, \Ac_{G^+ (\mathbb{Q})}  \,$
is given by the stated formula.
\end{proof}
\bigskip

\noindent \noindent We shall next analyze the above derivation, it
is given by construction by the 1-cocycle on $G^+ (\mathbb{Q})$
with values in $\Mc_4$,
\begin{equation} \label{sig}
\sigma_{\g}=X(\mu_{\g}) - \frac{1}{2}\, \mu_{\g}^2 \, .
\end{equation}
We recall that the classical Schwarzian, for which we use the
standard notation  $\{y \, ;x\}$, is given by the expression
\begin{equation} \label{Schw}
\{y \, ;x\} \, := \, \frac{d^2}{dx^2} \left( \log \frac{dy}{dx}
\right) \, - \, \frac{1}{2} \left( \frac{d}{dx} \left( \log
\frac{dy}{dx} \right) \right)^2 \, .
\end{equation}
The Schwarzian vanishes precisely on the M\"obius transformations
\begin{equation} \label{mob}
 z \mapsto \g \cdot z = \frac{az+b}{cz+d} \, ,
 \quad \g = \begin{pmatrix} a &b\\
 c &d  \end{pmatrix} \in SL(2,\mathbb{R})\, ,
\end{equation}
which in particular act on $S^1$ viewed as ${\mathbb P}^1
({\mathbb R})$. In the present context of modular forms we shall
show that the Schwarzian is an inner derivation of $\, \Ac_{G^+
(\mathbb{Q})}  \,$.
\medskip

\noindent The Schwarzian version of the chain rule can be
expressed as the identity
\begin{equation} \label{chain}
\{y \, ;x\} \, = \, \left( \frac{dz}{dx} \right)^2 \left(\{y \,
;z\} \, - \, \{x \, ;z\} \right)  \, ,
\end{equation}
where $ \, y = y(x) \,$ and $\, x = x(z) $. Using (\ref{mob}), it
follows that if
$$ Y (z) \, = \, ( y|_{0} \, {\g} ) (z) \,
= \, y(\frac{az+b}{cz+d}) \, ,
$$
then
\begin{equation} \label{Y}
\{Y \, ;z\} \, = \, \{y \, ;z\}|_{4} \, {\g} \, .
\end{equation}

\noindent Now, in view of (\ref{Xmu}), one has
\begin{equation} \label{X2}
\left( X(\mu_{\g})- \frac{1}{2}\, \mu_{\g}^2 \right) \cdot (2  \pi
i \, dz)^2 \, = \, \{Z |_{0} {\g}\, ;Z\}(dZ)^2 \, .
\end{equation}
By the above chain rule,
\begin{equation} \label{Z1}
\{Z |_{0} {\g}\, ; Z\}=(\frac{dz}{dZ})^2\ \;(\{Z |_{0} {\g} \,
;z\} - \{Z \, ;z\}) \, ,
\end{equation}
while by (\ref{Y}),
\begin{equation} \label{Z2}
\{Z |_{0} {\g} \, ; z \}=\{Z \, ; z\}|_{4} \, {\g} \, .
\end{equation}
Thus, from the identities (\ref{X2}) -- (\ref{Z2}) one obtains for
the ``Schwarzian'' 1-cocycle $\s$ defined by (\ref{sig}) the
formula
\begin{equation}
\sigma_{\g}=(2  \pi i)^{-2}\,(\{Z\, ;z\}|_{4} \, {\g}-\{Z \, ;z\})
\, .
\end{equation}
This shows that $\s$ is the coboundary of the weight $4$ modular
form
\begin{equation}
\om_4 \, = \, (2  \pi i)^{-2}\, \{Z \, ;z\} \, .
\end{equation}
A direct computation, using the definitions (\ref{Schw}) and
(\ref{ram}) together with (\ref{newZ}) gives
\begin{equation} \label{xg}
\om_4 \, = \, X (g_2^*  ) + \frac{1}{2} \, (g_2^*)^2 \, ,
\end{equation}
where
\begin{equation} \label{lg2}
g_2^* \, := \, \frac{1}{2  \pi i} \, \frac{d}{dz} (\log \eta^4 )
 \, = \, - \frac{2}{(2  \pi i)^2} \, G_2^* \, .
\end{equation}
In turn this implies (cf. e. g.~\cite{Za}) that $\om_4$ is a
multiple of the classical Eisenstein series $G_4$, namely
\begin{equation}
\om_4 \, = \, -\frac{E_4}{72} \, = \, - \frac{10}{(2  \pi
)^4}\,G_4 \, .
\end{equation}
where
\begin{equation}
E_4(q):= 1 + 240\,\sum_1^{\infty} n^3 \frac{q^n}{1-q^n}
\end{equation}

\noindent We can thus conclude with the following statement.
\medskip

\begin{proposition} \label{inner}
The derivation $\d'_2$ of the algebra $\, \Ac_{G^+ (\mathbb{Q})}
\,$ is inner and implemented by the weight $4$ modular form $\om_4
\, = - \frac{E_4}{72}$,
\begin{equation*}
\d'_2 (a) = [a, \, \om_4]\, , \quad \fl \, a \in \Ac_{G^+
(\mathbb{Q})}  \, .
\end{equation*}
\end{proposition}
\bigskip
\noindent The obvious question then, is whether one can perturb
the action of $\Hc_1$ on $\Ac$ by a $1$-cocycle in such a way
that, for the perturbed action, $\delta'_2=0$. The new action of
$\Hc_1$ on  $\Ac$ would then fulfill, more generally,
$$
\delta_n \, = \, \frac{(n-1)!}{2^{n-1}} \, \delta_1^n
$$
and therefore would come from a much smaller Hopf algebra.
\smallskip

\noindent We shall show however that such a perturbation does not
exist, and we shall actually relate the class of the action of
$\Hc_1$ modulo such perturbations to  the  data that
Zagier~\cite{Za} introduced to define canonical Rankin-Cohen
algebras.
\bigskip

\noindent To fix the notation, let us recall a few simple
generalities on cocycle perturbations of Hopf actions (see e.g.
\cite{miriam}).

\noindent  Let $\Ac$ be an algebra and $\Hc$ a Hopf algebra. We
let $\displaystyle \Lc( \Hc ,\Ac)$ be the convolution algebra of
linear maps from $\Hc$ to $\Ac$. The product of two elements of
$\displaystyle \Lc( \Hc ,\Ac)$ is given by
$$
u\, v \, (h) \, = \, \sum  u(h_{(1)})\; v( h_{(2)}) \,,\quad
\forall h \in \Hc \, ,
$$
using the standard short-hand notation, where one denotes the
coproduct by
$$
 \Delta(h) = \sum h_{(1)} \ot h_{(2)} \,,\quad \forall h \in \Hc \, .
$$
Let us also assume that the Hopf algebra $\Hc$ is acting, by a
Hopf action, on $\Ac$. A $1$-\textit{cocycle} is then an
invertible element of the convolution algebra of linear maps from
$\Hc$ to $\Ac$,
 $ \, u \in \Lc( \Hc ,\Ac) \, $,
such that
\begin{equation} \label{1coc}
u(h \, h') = \sum u(h_{(1)}) \, h_{(2)}(u(h')) \,,\quad \forall h
\in \Hc \, .
\end{equation}
Note that in view of the above equation, in the case of our Hopf
algebra, we only need to compute $u$ on the generators $X, Y,
\delta_1$ of $\Hc_{1}$ in order to determine it uniquely.
\smallskip

\noindent The perturbed action of $\Hc$ on $\Ac$ is given by,
$$
\tilde{h}(a):= \sum u(h_{(1)}) \, h_{(2)} (a) \, u^{-1}(h_{(3)})
$$
We look for a cocycle $u$ such that  $\tilde{\delta'_2}=0$ for the
perturbed action, or equivalently
$$
u(\delta'_2) \, = \, \om_4 \, .
$$
Since the action of $\Hc_1 $ under consideration actually commutes
with the natural coaction of $G^+ (\mathbb{Q})$ on $\Ac_{G^+
(\mathbb{Q})}$, it is natural to only consider cocycles $u$  which
preserve this property. It then follows that the values taken by
such a $1$-cocycle $u$ on generators must belong to the subalgebra
$\Mc \subset \Ac_{G^+ (\mathbb{Q})}$ and  have to be of the
following form:
$$
u(X)= t \in \Mc_2, \;u(Y)= \lambda \in \mathbb{C}, \; u(\delta_1)=
m \in \Mc_2 \, .
$$

\medskip

\begin{proposition} Let $t \in \Mc_2, \; \lambda \in \mathbb{C}, \;
 m \in \Mc_2 \,$.
\begin{itemize}
\item[$1^{0}$.] There exists a unique $1$-cocycle
$\displaystyle u \in \Lc( \Hc_1 ,\Ac_{G^+ (\mathbb{Q})})$ such
that
$$
u(X)= t , \;u(Y)= \lambda , \; u(\delta_1)= m  \, .
$$
\item[$2^{0}$.] The conjugate under $u$ of the
       action of $\Hc_1$ is given
       on generators as follows:
\begin{eqnarray}
\tilde{Y} = Y \, , \quad \tilde{X}(a) &=& X(a)+[(t- \lambda m),a]
-
 \lambda\, \delta_1(a)+ m Y(a) \, , \nonumber \\
 \tilde{\delta_1}(a) &=& \,\delta_1(a)
+ [m,a] \, , \qquad a \in \Ac_{G^+ (\mathbb{Q})}\,. \nonumber
\end{eqnarray}
\item[$3^{0}$.] The conjugate under $u$  of $\delta'_2$ is given
by the operator
$$
\tilde{\delta'_2}(a) = [X(m) +\frac{ m^2}{2}- \, \om_4,\; a] \, ,
\qquad a \in \Ac_{G^+ (\mathbb{Q})} $$ and there is no choice of
$u$ for which $\tilde{\delta'_2} = 0$.
\end{itemize}
\end{proposition}

\medskip
\begin{proof} $1^0$ The uniqueness is obvious but one needs to check
the existence. Any element of $\Hc_1$ can be uniquely written as a
linear combination of monomials of the form $\displaystyle
P(\d_1,\d_2,...,\d_k)X^n Y^l$. Let $m^{(k)}$ be defined by
induction by $\displaystyle m^{(1)}=m$ and
\begin{equation}
m^{(k+1)}=\, X(m^{(k)}) + m\, Y(m^{(k)})
\end{equation}
Similarly, let $t^{(k)}$ be defined by induction by $\displaystyle
t^{(1)}=m$ and
\begin{equation}
t^{(k+1)}=\, X(t^{(k)}) +\, m\, Y(t^{(k)})+\, t \; t^{(k)}
\end{equation}
Let us then define $u$ by the equality,
\begin{equation} \label{coc11}
u( P(\d_1,\d_2,...,\d_k)X^n Y^l):=\, \lambda^l \,
P(m^{(1)},m^{(2)},...,m^{(k)})\, t^{(n)}
\end{equation}
One checks directly that $u$ is a $1$-cocycle and that $u^{-1}$ is
given by
\begin{equation}
u^{-1}( P(\d_1,\d_2,...,\d_k)X^n Y^l): =\,(-\lambda)^l\,
P(n^{(1)},n^{(2)},...,n^{(k)})\, s^{(n)}
\end{equation}
where  $\displaystyle n^{(k)}:=-\, X^k(m)$ and $s^{(n)}$ is
defined by induction using $\displaystyle s^{(1)}:=-t + \,
\lambda \, m$ and
\begin{equation}
s^{(n+1)}:=\, X(s^{(n)}) +s^{(n)}(-t + \,    \lambda \, m) \,.
\end{equation}

\noindent $2^{0}$ One has as above
 $u^{-1}(X)= -t + \,    \lambda \, m$.
Since
$$
\Delta^{(2)}(X)= X \ot 1 \ot 1 +1 \ot X \ot 1 + 1 \ot 1 \ot X +
\delta_1 \ot Y \ot 1 + \delta_1 \ot 1 \ot Y + 1 \ot  \delta_1 \ot
Y
$$
it follows that  $\tilde{X}$ is given by
$$
\tilde{X}(a)=X(a)+[(t- \lambda m),a] - \lambda \delta_1(a) + m
Y(a)\,.
$$
The computation of
 $\tilde{Y}$ and $\tilde{\delta_1}$
is  straightforward.

\noindent $3^{0}$ Using the above formula (\ref{coc11}) gives
\begin{equation} \label{*}
u(\delta'_2)=\, X(m) +\frac{ m^2}{2}
\end{equation}
We are thus reduced to showing that
 one cannot find a modular form (of arbitrary
level) which solves the equation
\begin{equation} \label{**}
X(m) +\frac{ m^2}{2}= \, \om_4 \, .
\end{equation}
To do that let us use the ``quasimodular'' (\cite{Zagier})
   solution
to the above equation, given by $g_2^*$ (\ref{lg2}). Note that not
only $g_2^*$ fulfills (\ref{*}) but also one has
$$
X(f) =\frac{1}{2  \pi i} \frac{d}{dz}\, f -g_2^* \, f \, ,
$$
for any modular form of weight 2. Thus,
$$
\frac{1}{2  \pi i} \frac{d}{dz}\,g_2^*- \frac{1}{2} (g_2^*)^2 = \,
\om_4
$$
so that one gets a minus sign and not a plus sign when trading $X$
for $ \displaystyle q \, \frac{d}{dq}$.

\noindent Assuming now that $f$ is a solution of (\ref{**}), one
gets
$$
\frac{1}{2  \pi i} \frac{d}{dz}\,f -\,g_2^* f + \frac{1}{2} f^2 =
\frac{1}{2 \pi i} \frac{d}{dz}\,g_2^*- \frac{1}{2}( g_2^*)^2 \, .
$$
In turn, this implies
$$
\frac{1}{2 \pi i} \frac{d}{dz}(f-g_2^*) +\frac{1}{2} (f-g_2^*)^2=0
\, .
$$
But one can trivially integrate this non-linear equation. Indeed,
the substitution $\displaystyle y = \frac{1}{f - g_2^*}$ gives
$$
\frac{1}{2 \pi i} \frac{d}{dz}\, y \, = \, \frac{1}{2} \, ,
 \qquad \text{so that} \qquad
f- g_2^* \, = \, \frac{1}{ \pi i z + c} \, \, ;
$$
the latter fails to be periodic in $z$ except when it vanishes.
\end{proof}

\bigskip
\noindent \textbf{Remark 2. \,} Note  that the freedom we have in
modifying the action of the Hopf algebra $\Hc_1$ by a $1$-cocycle
 exactly modifies the restriction of
the action of $X$ on modular forms
 and the value of $\om_4$ as in the
data used by Zagier~\cite{Za} to define ``canonical'' Rankin-Cohen
algebras. More precisely the derivation $\partial$ and the element
$\Phi$ of degree $4$ used by Zagier in \cite[Proposition 1]{Za}
 correspond
to the restriction of $X$ to $\Mc$ and to $\om_4 \, = \, 2\Phi$.
The gauge transformation associated in \textit{loc. cit.} to an
element $\phi$ of degree 2 corresponds to the value $\,
\displaystyle m  =  - \frac{\phi}{2} \,$ for our cocycle.
\bigskip

\noindent A first step in the understanding of the relationship
between actions of the Hopf algebra $\Hc_1$ such that the
derivation $\d'_2$ is inner and `canonical' Rankin-Cohen algebras
consists in statement $3^{0}$ of the following theorem, which we
shall formulate for the modular Hecke algebra $\Ac (\G)$
associated to an arbitrary congruence subgroup $\G$ of $\SL (2,
\mathbb{Z})$.
\smallskip

\noindent We shall in fact show that all the above results also
apply to the algebras $\Ac (\G)$. In view of the adelic
interpretation of $\Ac (\G)$ as a reduced algebra of the crossed
product of $\Mc$ by  $  \GL (2 , \mathbb{A}_{f})^0 $, the proofs
will be similar, provided one checks that the reducing projection
$e_{\G}$ is invariant under the action of $\Hc_1$. It is more
instructive however to proceed directly and concretely describe
the action of $\Hc_1$ on $\Ac (\G)$, as follows below.
\smallskip

\noindent The action of the generators of $\Hc_1$ on $\Ac (\G)$ is
given by the following operators:
\begin{eqnarray}
 Y(F)_{\a} \,&:=& \, Y(F_{\a}) \, , \qquad \fl \,
 F \in \Ac (\G) \, , \a \in G^+ (\mathbb{Q}) \, ,  \cr
  X(F)_{\a} \, &:=& \, X(F_{\a}) \, ,  \qquad   \cr
 \d_n (F)_{\a}  &:=& \, \mu_{n, \a} \cdot F_{\a} \, , \label{actn}
\end{eqnarray}
where
\begin{equation}
\mu_{n, \, \a} \, := \, X^{n-1}(\mu_{\a}) \, ,
 \quad \fl n \in {\mathbb N} \, .
\end{equation}
The right hand side of each of the formulae in (\ref{actn}) gives
an element of $\Ac (\G)$. This is true for the second  because
$\mu_{\g}$ vanishes for $\g \in \G \sbs \G(1)$, and it follows for
the third since, cf. (\ref{zmu}),
\begin{equation*}
 \mu_{n, \, \a \g} \, = \mu_{n, \, \a} \, , \quad \g \in \G \sbs \G(1) \, .
\end{equation*}
\noindent The  algebra $\Ac (\G)$ is graded by the weight of
modular forms, i. e. by the operator $Y$. The weight $0$
subalgebra given by $Y=0$, is  the standard Hecke algebra
associated to  $\G$.
\medskip

\medskip

\begin{theorem} \label{trans}
    Let $\G$ be any congruence subgroup.
\begin{itemize}
\item[$1^0$.] The formulae (\ref{actn}) define a Hopf action of
the Hopf algebra $\Hc_1$ on the algebra $\Ac (\G)$.

\item[$2^0$.] The Schwarzian derivation $\d'_2=\d_2 -\frac{1}{2}\d_1^2$
 is inner and is
implemented by $\, \om_4 \in \Ac (\G) \, $.

\item[$3^0$.] The image
of the generator (\ref{fund}) of the tranverse fundamental class
 $\, [F] \in HC_{\Hopf}^{2}(\Hc_1) \,$
under the canonical map from the Hopf cyclic cohomology of $\Hc_1$
to the Hochschild cohomology of $\Ac (\G)$,
\begin{equation*}
 \chi (F) (a_1, a_2)  :=  X(a_1) Y(a_2) - Y(a_1) X(a_2) -
 \d_1 (Y(a_1)) Y(a_2) \, ,
\end{equation*}
gives the natural extension of the first Rankin-Cohen bracket $\{
\, \cdot\,\, ,\, \cdot \}_1$ to the  algebra  $\Ac (\G)$.

\item[$4^0$.] The cuspidal ideal $\Ac^0 (\G)$ is globally invariant
under the action of $\Hc_1$.
\end{itemize}
\end{theorem}
\medskip

\begin{proof} The verification of the first two statements
is identical to that already given for $\Ac_{G^+ (\mathbb{Q})}$.
The third statement follows from $1^0$ and the cocycle property of
$F$, together with the formula
$$
\{ \, a_1\, ,\, a_2 \}_1=X(a_1) Y(a_2) - Y(a_1) X(a_2)
$$
for the Rankin-Cohen bracket of modular forms. The fourth
statement is true because $X$ preserves the cuspidal ideal
$\Mc^0$.
\end{proof}
\bigskip

\noindent  Note that, exactly as before, inner perturbations of
the action of $\Hc_1$ by a cocycle of the form
$$
u(X)= t \in \Mc_2, \;u(Y)= \lambda \in \mathbb{C}, \; u(\delta_1)=
m \in \Mc_2 \, .
$$
yield a new action of $\Hc_1$ with
\begin{equation*}
\tilde{Y} \, = \, Y \, , \quad \tilde{X}\, = \, X + \ad (t -
\lambda m) -
 \lambda \delta_1 + m Y \, , \quad
 \tilde{\delta_1} \, = \,\delta_1
+ \ad (m ) \, ,
\end{equation*}
and with the new Schwarzian $\tilde{\delta'_2}$ inner, implemented
by $\, \displaystyle \om_4 - X(m) - \frac{ m^2}{2} \,$.

\noindent In the simplest case $\G= \G(1)$ one can witness the
non-triviality of the action of $\Hc_1$ on $\Ac (\G)=\Ac (1) $ by
computing the higher derivatives $\d_k(T_n)$ of the standard Hecke
operators $T_n$.
\bigskip

\noindent We conclude this section by describing in precise terms
the analogy between the action of $\Hc_1$ on crossed products of
polynomial functions on the frame bundle of $S^{1}$ by discrete
subgroups of $\Diff (S^{1})$ (cf. \cite{CM1}) and its action on
modular Hecke algebras.
\medskip

\noindent Since $\eta^4$ has level $6$, the natural domain of
definition of the differential form defined by equation
(\ref{newZ}) is the modular elliptic curve
$$ E \, := \, X(6) \, \cong X_0 (36) \, ,
$$
where the isomorphism comes from the fact that $\G(6)$ is
conjugate to $\G_0 (36)$ by the matrix
$ \displaystyle  \begin{pmatrix} 6 &0\\
0 &1  \end{pmatrix} \in G^+ (\mathbb{Q}) $. The operator $X$ of
equation (\ref{ram}) is simply giving the canonical flat
connection on the canonical line bundle of $E$ .
Proposition~\ref{gend} then shows that the action of $\Hc_1$ on
$\Ac$ is governed by exactly the same rules as its action on
crossed products for the case of $\Diff (S^1)$ (cf. \cite{CM1}).
By Theorem~\ref{trans} the same holds for the action of $\Hc_1$ on
the modular Hecke algebras $\Ac (\G)$, where the diffeomorphisms
are replaced by Hecke correspondences.  Thus, both `differential'
operators $X$ and $Y$ have clear geometric meaning, while the
$\d_n \,$, $\, n \in \mathbb{N}$, are uniquely determined by the
Hopf action rules.
\smallskip

\noindent  Pursuing the above analogy we shall now describe the
projective structure on the curve $E$ which is responsible for the
fact that the derivation $\delta'_2$ of the modular Hecke
algebras (Theorem~\ref{trans}) is inner. The simplest projective  structure
on $S^1$ is given by the familiar identification $S^1 \cong
\mathbb{P}^1 (\mathbb{R})$. The relation between the angular
variable $\theta \in \mathbb{R}/2 \, \pi \,\mathbb{Z}$ and the
variable $t \in \mathbb{R} $ of the rational parametrization
defining the  projective structure is $\displaystyle \, t \, = \,
\tan (\frac{\theta}{2}) \,$, which shows that, when expressed in
terms of $\theta$, this projective structure on  $S^1 $ is given
by the quadratic differential
\begin{equation*}
    \rho \, : = \, \{t\, ; \theta \} \, d \theta^2 \, = \,
    \frac{1}{2} \, d \theta^2 \, .
\end{equation*}
In our case the elliptic curve $E$ plays the role of the circle
$S^1$, $Z$ plays the role of the angular variable $\theta$ and the
projective structure is defined by the variable $z$, and hence
given by the quadratic differential
\begin{equation} \label{qd}
    \varpi \, : = \, \{z\, ; Z \} \, d Z^2 \, = \,
    - \{ Z\, ; z \} \, d z^2 \, .
\end{equation}
Now the elliptic curve $E$ has Weierstrass equation
\begin{equation} \label{W}
    y^2 \, = \,  x^3 + 1
\end{equation}
and in these terms,
\begin{equation} \label{Zx}
    dZ \, = \, \frac{dx}{2\,y}
\end{equation}
We need to express the above quadratic differential in terms of
its Weierstrass coordinates $x$ and $y$. Since by equation
(\ref{xg})
\begin{equation*}
     \{ Z\, ; z \} \, = \, ( 2 \pi i)^2 \, \om_4
   \, = \, - \frac{( 2 \pi i)^2 }{72} \, E_4 \, ,
\end{equation*}
what is needed is to express the ratio $\displaystyle \frac{E_4}{
\eta^8 } \,$ as a rational function of  $x$ and $y$. The function
$x$ is the only solution of the differential equation
\begin{equation} \label{xxx}
    ( \frac{dx}{dZ})^2  \, = \, 4\,( x^3 + 1)
\end{equation}
having an expansion near $q=0$ of the form
\begin{equation} \label{dx}
    x \, = \, q^{-\frac{1}{3}} \, (1+ \sum_{n = 1}^{\infty} a_{n}
    q^{n} ) \, .
\end{equation} \label{x}
It is given by the ratio
\begin{equation*}
    x \, = \, \frac{\mu}{\eta^8} \,  ,
\end{equation*}
where  $\, \mu \,$ is the  $\G_0 (6)$-modular form of weight $4$
\begin{equation*}
    \mu \, = \, \frac{1}{5} \, \left( g_{4}
    -  4 \,  g_{4} \vert \begin{pmatrix} 2 &0 \\ 0 &1  \end{pmatrix}
     -  9 \,  g_{4} \vert \begin{pmatrix} 3 &0 \\ 0 &1  \end{pmatrix}
      +  36 \, g_{4} \vert \begin{pmatrix} 6 &0 \\ 0 &1  \end{pmatrix}
      - 36 \, S_{\rm new} \right) \,
\end{equation*}
with $ g_{4}:= \frac{1}{240}\,E_4\,$ and
\begin{equation*}
     S_{\rm new} (q) \, = \, q - 2q^2 - 3q^3 + 4q^4 +6q^5 +6q^6
    - 16q^7 + \ldots
\end{equation*}
denoting the new form of weight $4$ and level $\G_0 (6)$.
\smallskip

\noindent The transformation $\, \a = \begin{pmatrix} 1 &3 \\ 0 &1
\end{pmatrix} \in \G_0 (6) \,$ defines an involution on our curve
$E$, with respect to which $\, \eta^4 \vert \a = - \eta^4 $. Thus,
with $\displaystyle 2\, y \, := \, \frac{dx}{dZ} \,$, one has
$$ x \vert \a = x \, , \text{and} \quad  y \vert \a = - y \, .
$$
On the other hand,
$$  E_4 \vert \a \, = \, E_4 \, ,
$$
which shows that $\displaystyle \frac{E_4}{ \eta^8 } \,$ is a
rational function of the above $x$. Explicitly, this rational
function $R(x)$ can be computed as
\begin{equation} \label{ratfrac}
     \frac{E_4}{ \eta^8 }  \, = \,R(x) = \,
 \frac{ (x^3+4)(x^9 +228\, x^6+48 \,x^3+64) }{ x^2 ( x^3 - 8)^2 (x^3 +1)} \, ,
\end{equation}
where the poles correspond precisely to the $12$ cusps of $ E =
X(6)$.
\smallskip

\noindent From the equation (\ref{newZ}), one has
\begin{equation*}
  \varpi \, = \,   \frac{(2\pi i)^2}{72} \, E_4 \, (dz)^2
  \, = \, \frac{1}{2} \, \frac{E_4}{ \eta^8 }  \, (dZ)^2
\end{equation*}
and using equation (\ref{xxx}) in the form
\begin{equation*}
 4\,(dZ)^2 \, = \,   \frac{(dx)^2}{x^3 + 1} \,
\end{equation*}
we finally obtain
\begin{equation} \label{pr}
  \varpi \, = \,
 \frac{(x^3+4)(x^9 +228\, x^6+48 \,x^3+64)}{8 \,(x(x^3-8)(x^3+1))^2}
\, dx^2  \, .
\end{equation}

\smallskip

$$
\hbox{ \psfig{figure=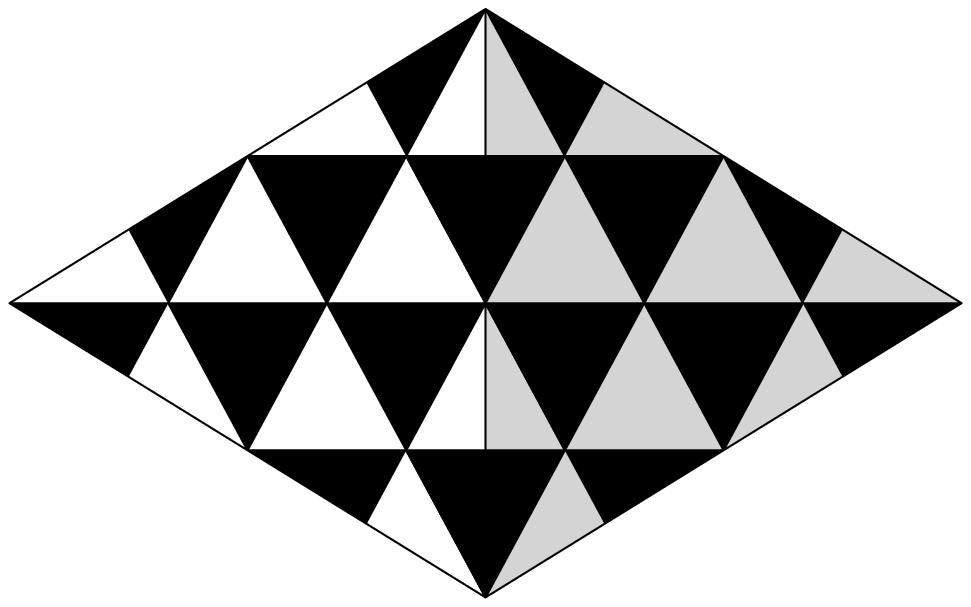} }
$$

\centerline{The Elliptic Curve E and the Subgroup $C$}

\bigskip \bigskip
\noindent The curve $E$ is the simplest elliptic curve in that it
admits a triangulation by two equilateral triangles (the two big
triangles in the picture). It is well-known, cf.~\cite[Theorem
I.5.2]{Bost}, that the existence of an equilateral triangulation
defining the conformal structure characterizes curves defined over
$\bar{\mathbb{Q}}$ among compact connected Riemann surfaces. The
group  $\PSL(2,\mathbb{Z} /6)$ acts by automorphisms of $E$ and its
commutator subgroup $C \cong (\mathbb{Z} /3) \cdot (\mathbb{Z}
/2)^2$ preserves the $1$-form  $\eta^4 dz$ and hence acts by
translations, and can thus be identified with a subgroup of the
abelian group $E$. The corresponding elements of $E$ are the 12
cusps. They provide the vertices of  a triangulation of $E$ by
$24$ equilateral triangles which is invariant by the group $C$ and
is displayed in the above picture. The quotient $E/C$ is the
elliptic curve $\G'(1) \bash H^* $ i.e. the quotient  $\mathbb{C}
/ \Lba \,$ which is again triangulated by two equilateral
triangles. One then gets  the following geometric expression for
the rational fraction $R(x)$ of (\ref{ratfrac}),
\begin{equation}
R(x)=\sum_{\tau \in C} x^{\tau }
\end{equation}
where $x^{\tau }$ is the transform of $x$ by the translation
$\tau $ of the elliptic curve. Of course $\eta^4 dz$ does make
sense on  the quotient $E/C = E'$ and the formula for the
projective structure looks simpler if we work directly there; it
becomes
\begin{equation} \label{pr2}
  \varpi' \, = \, \frac{x \,  dx^2 }{2 \, y^2} \, = \,
 \frac{x \,  dx^2 }{8 \,(x^3-1728)} \, .
\end{equation}

\noindent In order to recast the above developments in a
completely geometric perspective one would still need to be able
to characterize the Hecke correspondences $\sigma(\alpha)$ over
$E$ associated to
 $\alpha \in \G\bash G^+ (\mathbb{Q})/\G \,$
(cf. \cite[7.2.3]{Sh}) in geometric terms. Letting $\Lc$ be the
ample   line bundle associated to the divisor of cusps on $E$, the
triple
$$
(E, \sigma(\alpha), \Lc)
$$
is then strikingly similar to the geometric data associated to
quadratic algebras (see \cite{ATV} for instance) with the
important nuance that instead of the graph of  a translation on
the elliptic curve we are dealing in our context with the
irreducible curve representing a Hecke correspondence.

\bigskip

\section{The transverse cocycle and the Euler class}

\noindent In this section we shall concentrate on the analogue of
the Godbillon-Vey cocycle obtained from the cyclic cocycle $\,
[\d_{1}] \, \in \, HC_{\Hopf}^{1} \, (\Hc_1) \,  $
 (cf. Proposition \ref{gv2})
 in the above
action of $\Hc_1$ on the crossed product algebras. We need for
that a substitute for the invariant trace $\tau$ which was used in
the context of codimension $1$  foliations. We first  consider the
projection map (already used in the \textit{residue} definition in
\cite{CMZ})
$$ P : \Mc \, \ra \, \Mc_{2} \, ,
$$
which projects onto the component of weight 2, and promote it to a
projection denoted by the same letter $\, P : \Ac_{G^+
(\mathbb{Q})}  \, \ra \, \Mc_{2} \,$, by setting
\begin{equation*}
P (f \, U_{\g}^*) \, = \,\begin{cases}
 P(f)  &\text{ if }  \g = 1  \\
0 &\text{ if }  \g \ne 1
\end{cases}
\end{equation*}
The next result expresses its covariance under the action of
$\Hc_1$.
\bigskip

\begin{lemma} \label{cov}
For any $\, a ,  b \, \in  \Ac_{G^+ (\mathbb{Q})}  \,$ and any $\,
h \in \Hc_{1} \,$, one has
\begin{equation} \label{inv}
P \, (h(a) \cdot b) = P \, (a \cdot \wt S (h) \, (b)) \, ,
\end{equation}
where $\wt S = \nu * S$ is the twisted antipode.
\end{lemma}
\smallskip

\begin{proof} One can assume that $a$, $b$ are of the
form $a=f \, U_{\g}^*$, $b=g \, U_{\g^{-1}}^*$ and that  $h$ is a
monomial in the generators of $\Hc_{1}$. We endow $\Hc_{1}$ with
its natural grading so that deg($Y$)$=0$, deg($X$)$=1$,
deg($\d_n$)$=n$. Since $P(f) = 0$ unless the weight of $f$ is $2$,
both sides of (\ref{inv}) vanish unless
$$ {\rm weight}(f) + {\rm weight}(g) + 2\, {\rm deg}(h)=2 \, .
$$
If ${\rm deg}(h)\neq 0$, both $f$ and $g$ are of weight $0$ and
hence are constants. Then both sides of (\ref{inv}) vanish unless
the monomial is of the form $Y^m \d_1$. In that case, we have
$h(a)=\mu_{\g}\, U_{\g}^*$ and $h(a) \cdot b=\mu_{\g}$. Since
$$\, \wt{S} (Y) = -Y + 1 \, \,  \text{and}
 \, \, \wt{S} (\d_1) = -\d_1 \, ,
 $$
we get
$$ \wt{S} (Y^m  \d_1) = -\d_1 (1-Y)^m \, \,
\text{and} \, \, \wt{S} (h)(b)=-\mu_{\g^{-1}}U_{\g^{-1}}^* \, .
$$
Then
$$ a \cdot \wt{S} (h)(b)=-\mu_{\g^{-1}}|_{2}  \g
$$
and thus the cocycle property of $\mu$ gives (\ref{inv}) for the
case considered.

\noindent If ${\rm deg}(h)= 0$, then $h$ is of the form $h= Y^m$.
Let then $ \, f \in \Mc_k \,$. Both sides of (\ref{inv}) vanish
unless $ \, g \in \Mc_{2-k} \,$. Since $\, \wt{S} (Y) = -Y + 1 \,
$, it follows that
$$ P (Y^m(a) \,b) = \, \left(\frac{k}{2} \right)^m a \, b \, ,
$$
while
$$  P (a \cdot \wt{S} (Y^m)(b))= \,
\left( 1 - \frac{2-k}{2} \right)^m  a \, b \, ,
$$
which gives the required equality.
\end{proof}
\bigskip

\noindent At this point, inspired by the construction of the
Godbillon-Vey cocycle, we set
$$
\GV \, (a, b) \, = \, P (a \cdot \d_{1} (b) ) \, , \qquad \fl \,
a,b \in \Ac_{G^+ (\mathbb{Q})}  \, .
$$
In the absence of an invariant trace, the construction of a formal
analogue of the Godbillon-Vey cocycle will be completed by means
of coupling with the analogue of the classical period cocycle. Due
to the specificity of the situation, the entire construction can
be conveniently phrased in terms of group cohomology rather than
cyclic cohomology. Thus, $\, \GV \,$ reduces to the group
1-cocycle $ \, {\Eb} \in Z^1 (G^+ (\mathbb{Q}), \Mc_2) $,
\begin{equation}
{\Eb} (\g) \, : = \, \GV \, ( U_{\g^{-1}}^* , U_{\g}^* ) \, = \,
\mu_{\g}|_{2}  \g^{-1}   \, , \qquad \fl \, \g \in G^+
(\mathbb{Q}) \, ,
\end{equation}
which we call the {\it transverse  cocycle}. This cocycle $\,
[{\Eb}] \in H^{1} \, ( G^+ (\mathbb{Q}), \Mc_2 ) \,$ is
nontrivial. We can afford to omit the easy and direct proof, since
the subsequent results will yield a finer understanding.
\bigskip

\noindent In order to provide a conceptual framework for the next
proposition we shall now describe in details the identification of
the $1$-jet bundle  $J^1 (H_\mathbb{A})$ with the line bundle
$\Lc_\mathbb{A}^{-2}$ where $\Lc_\mathbb{A}$ is the lattice line
bundle.
\smallskip

\noindent A modular form of weight $2k$ is a differential form
$f(z)\, dz^k$ on the complex curve $H_\mathbb{A}$, and hence can
be viewed directly  as a polynomial function on the (complex)
1-jet bundle $J^1 (H_\mathbb{A})$. The careful description of the
isomorphism $J^1 (H_\mathbb{A}) \sim \Lc_\mathbb{A}^{-2}$ does
give a very useful description of the invariant  connections  on
these line bundles. We use the canonical identification of the
Poincar\'e disk with the space of translation invariant conformal
structures on $\mathbb{C}$, by means of Beltrami differentials.
Thus, with $\mu$ in the Poincar\'e disk, the conformal structure
on $\mathbb{C}_{\mu}$ is defined by declaring that $dz + \mu
d\bar{z} $ is a 1-form of type (1,0). The 1-forms $\{ d\mu, dz
+ \mu d\bar{z} \}$ are then a basis of forms of type (1,0) on the
total space of the pairs $(z, \mu)$. Viewing a fixed lattice $\,
\Lba \,$ as a lattice in each of the complex lines
$\mathbb{C}_{\mu}, \, \vert \mu \vert<1$, defines a map from
lattices to $1$-parameter families of lattices. Equivalently, this
amounts to associate to a given lattice $\, \Lba \,$ the following
germ of map from $\mathbb{C}$ to lattices $\Lba( \mu)$:
$$
\Lba( \mu):= \{ \om + \frac{\mu}{{\rm cov}(\Lba)}\, \bar{\om}
,\quad \om \in \Lba \} \, ,
$$
where ${\rm cov}(\Lba)$ is $-2i$ times the covolume of $\Lba$ and
$\mu$ is a complex number such that
$$
\vert \mu \vert< \vert{\rm cov}(L)\vert \, .
$$
The $\mathbb{Q}$ structure of $\Lba( \mu)$ is trivially obtained
from that of $\Lba$. The point of the normalization by the
covolume is that one has
$$
(a\, \Lba)(\mu)=a \,\Lba(\frac{\mu}{a^2}) \, , \qquad \forall \, a
\in \mathbb{C}^{\ts} \, ,
$$
which is exactly the correct homogeneity to identify the $-2$
power of the bundle of lattices with the 1-jet bundle.
\bigskip

\noindent For $\Lba = \mathbb{Z}\om_1 + \mathbb{Z} \, \om_2 $, one
has
$$
\Lba(\mu) = \mathbb{Z}(\om_1 + \frac{\mu}{{\rm cov}(\Lba)}\,
\bar{\om}_1) + \mathbb{Z} \, (\om_2 + \frac{\mu}{{\rm
cov}(\Lba)}\, \bar{\om}_2)
$$
and
$$
{\rm cov}(\Lba)= \om_2 \, \bar{\om}_1 -\om_1 \, \bar{\om}_2 \, .
$$
The inhomogeneous coordinate $z=\om_1/\om_2$ of $\Lba(\mu)$ is
thus given by
$$
z(\mu)=\frac{\om_1(\om_1 \, \bar{\om}_2- \om_2 \, \bar{\om}_1) -
\mu \, \bar{\om}_1 } {\om_2(\om_1 \, \bar{\om}_2- \om_2 \,
\bar{\om}_1) - \mu \, \bar{\om}_2} \, ,
$$
The Taylor expansion of $z(\mu)$ at $\mu=0$ gives,
$$
z(\mu)=\om_1/\om_2 + \,  \mu /\om_2^2 + \mu^2 \, \bar{\om}_2 /
(\om_2^3(\om_1 \, \bar{\om}_2- \om_2 \, \bar{\om}_1))+...
$$
where the terms of degree two and higher are not holomorphic in
the lattice $\Lba$. The 1-jet however is holomorphic in $\Lba$
which gives the required canonical identification
$$
J^1 (H_\mathbb{A})\sim \,\Lc_\mathbb{A}^{-2} \, .
$$
The above construction yields a canonical connection on these line
bundles. In terms of homogeneous lattice
 functions $F$ of weight $2k$, the connection
 is given by definition as follows:
$$
\nabla_{\mathbb{R}}(F) (\Lba):=  \, \frac{1}{2\pi i}
 \frac{d}{d\mu} (F(\Lba(\mu)))\vert_ {\mu=0} \, .
$$
Note that $ \nabla_{\mathbb{R}}(F)$ is homogeneous of weight
$2k+2$ but is not in general a holomorphic function of $\Lba$,
even when $F$ is holomorphic. In the inhomogeneous coordinate
$z=\om_1/\om_2$ one gets,
$$
 \nabla_{\mathbb{R}}=  \frac{1}{2\pi i}
\frac{d}{dz}+\frac{2k}{2 \pi i(z-\bar{z})} \, ,
$$
which one recognizes as a standard connection in the theory of
modular forms (\cite{Zagier}). It is the canonical connection on
this holomorphic line bundle, associated to its natural invariant
Hermitian metric. This connection has non-zero curvature, given by
the area element. By its very construction, the connection
$\nabla_{\mathbb{R}}$ is invariant under $G^+ (\mathbb{Q})$ and
its relation to the flat connection $\nabla$ is given by
\begin{equation} \label{differ}
X=\nabla_{\mathbb{R}}- 2\, \phi_{\0b}\,  Y \, .
\end{equation}
Here $\phi_{\0b}= \frac{1}{4 \pi^2}G_{\0b}$, with  $G_{\0b}$
denoting the modular -- but nonholomorphic -- weight $2$
Eisenstein series
$$
\, G_{\0b} (z) \, = \, G_{\0b} (z, 0) \, ,
$$
obtained by taking the value at $\, s = 0 \,$ of the analytic
continuation of the series
\begin{eqnarray}
G_{\0b} (z, s) \, &=& \sum_{(m, n) \in \mathbb{Z}^2 \setminus 0}
(m z + n)^{-2} \, |m z + n|^{-s}  \nonumber \\
  &=& 2 \zeta (2 + s) \, + \,
2 \sum_{m \geq 1} \, \sum_{n \in \mathbb{Z}} (m z + n)^{-2} \, |m
z + n|^{-s} \, , \quad \Re s > 0 \, . \nonumber
\end{eqnarray}
It is related to $\, G^*_2 \,$ by the identity
\begin{equation} \label{GG*}
G^*_2 (z) \, = \, G_{\0b} (z) \, + \, \frac{2  \pi i}{z - \bar{z}}
\, ,
\end{equation}
which makes the check of (\ref{differ}) straightforward. The $G^+
(\mathbb{Q})$-invariance of $\nabla_{\mathbb{R}}$ then entails the
following result.

\bigskip

\begin{proposition} \label{emu0}
For any $\, \g \in G^+ (\mathbb{Q})  \,$ one has
\begin{equation}  \label{gb0}
 \mu_{\g} \, = \,  2 ( \phi_{\0b} | \g \, - \, \phi_{\0b} )
\end{equation}
\end{proposition}
\smallskip

\noindent In order to obtain a better description of the right
hand side of (\ref{gb0}) we rely on Hecke's construction
in~\cite[\S 2]{He} of canonical generators for the
$\mathbb{Q}$-vector space $\, \Ec_2 (\mathbb{Q}) \,$ of rational
Eisenstein series. An Eisenstein series $\, E \in \Ec_2
(\mathbb{Q}) \,$ if the constant term of the $q$-expansion of $\,
E \,$ at each cusp is a rational number.
\smallskip

\noindent We shall closely follow the exposition in \cite[\S
2.4]{St}. For each $\, {\ab} = (a_1, a_2) \in ( \mathbb Q \slash
\mathbb Z )^2  $ let
\begin{equation*}
 G_{\ab} (z, s) := \, \sum_{ {\mb} \in {\mathbb Q}^2 \setminus {\0b},
\, {\mb} \equiv {\ab} \, \text{mod} \, {\mathbb Z} }
 (m_{1} z + m_{2})^{-2} |m_{1} z + m_{2}|^{-s} \, , \quad \Re s > 0 \, \, ;
\end{equation*}
for fixed $z \in H , \, G_{\ab} (z, s) $ may be analytically
continued to a meromorphic function in the complex $s$-plane which
is holomorphic at $ s = 0 $, so that one can define
\begin{equation*}
 G_{\ab} (z) := \,  G_{\ab} (z, 0) \, .
\end{equation*}
One then has
\begin{equation*}
  G_{\ab} | \g \, = \,  G_{{\ab} \cdot \g} \, , \quad \fl \, \g \in
  \G (1) \, ,
\end{equation*}
which shows that $ G_{\ab} (z) $ behaves like a weight $2$ modular
form of some level $N$. It fails to be holomorphic, but only in a
controlled way\, ; more precisely, the function
\begin{equation*}
    z \, \mpo  \, G_{\ab} (z) \, + \, \frac{2 \pi i}{z - \bar{z}}
\end{equation*}
is holomorphic on $H$. Moreover, for $\, {\ab} = (a_1, a_2) \in (
\mathbb Q \slash \mathbb Z )^2 \setminus {\0b} $,
\begin{equation*}
     \wp_{\ab} (z) \, = \,  G_{\ab} (z) \, - \,  G_{\0b} (z)
\end{equation*}
is the ${\ab}-$division value of the Weierstrass $\wp$-function
and the collection of functions
\begin{equation*}
 \left\{ \wp_{\ab} \, \, ; \,
\, {\ab} \in \left( \frac{1}{N} \mathbb Z \slash \mathbb Z
\right)^2 \setminus {\0b} \right\}
\end{equation*}
generates the space of weight $2$ Eisenstein series of level $N$.
\smallskip

\noindent In order to obtain a basis for $\Ec_2 (\mathbb Q)$, one
considers for each $\, {\xb} = (x_1, x_2) \in
 \left( \frac{1}{N} \mathbb Z \slash \mathbb Z \right)^2 $
the additive character $\, \chi_{\xb} :
 \left( \frac{1}{N} \mathbb Z \slash \mathbb Z \right)^2 \ra
 {\mathbb C}^{\times} $ defined by
\begin{equation*}
  \chi_{\xb} \left( \frac{\ab}{N} \right) \, := \,
 e^{2  \pi i\,  (a_2 x_1 - a_1 x_2 ) }
\end{equation*}
and one forms the `twisted' Eisenstein series
\begin{equation} \label{basis}
  \phi_{\xb} (z) \, := \, \sum_{{\ab} \in
 \left( \frac{1}{N} \mathbb Z \slash \mathbb Z \right)^2} \,
 \chi_{\xb} ({\ab}) \cdot  G_{\ab} (z) \, .
\end{equation}
 The definition is independent of $N$, and if $ {\xb} = (x_1, x_2) \in
 \left( \frac{1}{N} \mathbb Z \slash \mathbb Z \right)^2 \setminus
 {\0b} $ then $ \, \phi_{\xb} \, $ is a weight $2$ Eisenstein series
 of level $N \,$.
The function $\,  \phi_{\xb} \,$ is also known in the literature
as the logarithmic derivative of the generalized Dedekind function
$\, \eta_{\xb} \,$,
\begin{equation*}
  \phi_{\xb} \, = \,
 \frac{1}{2  \pi i} \cdot \frac{d}{dz} \log \eta_{\xb}  \, .
\end{equation*}
\smallskip

\noindent To account for the special case when $ \, {\xb} =  {\0b}
\,$, one adjoins the non-holomorphic but fully modular function $
\, \phi_{\0b} \,$ given as above by
\begin{equation} \label{b0}
2 \pi i \cdot \phi_{\0b} (z) \, - \, \frac{1}{z - \bar{z}} \, = \,
2 \cdot \frac{d}{dz} \log \eta \, .
\end{equation}
\medskip

\noindent All the linear relations among the functions $\,
\phi_{\xb} \, , \, {\xb} \in ( \mathbb Q \slash \mathbb Z )^2  \,$
are encoded in the {\it distribution property}
\begin{equation} \label{dis}
  \phi_{\xb} \, = \, \sum_{{\yb} \cdot\check{\g} = {\xb}} \,
  \phi_{\yb} | \g  \, .
\end{equation}
where
$$
        \check{\g} \, = \, \det \g \cdot {\g}^{-1} \, .
$$
This allows to equip the extended Eisenstein space
$$
\Ec^*_2 (\mathbb Q) \, = \, \Ec_2 (\mathbb Q) \, \oplus \, \mathbb
Q \cdot \phi_0 \, ,
$$
with a linear $\, \text{PGL}^+ (2, \mathbb{Q})$-action, as
follows. Denoting
$$  \Sc \, := \, ( \mathbb Q \slash \mathbb Z )^2 \, ,
\quad \text{resp.} \quad  \, {\Sc}' \, := \, \Sc \setminus {\0b}
\,
$$
and identifying in the obvious way
$$ \text{PGL}^+ (2, \mathbb{Q}) \, \cong
\, M_2^+ (\mathbb{Z}) \slash \{ \text{scalars} \} \, ,
$$
where $\, M_2^+ (\mathbb{Z})  \, $ stands for the set of integral
$2 \times 2$-matrices of determinant $ > 0$, one defines the
action of $\, \g \in M_2^+ (\mathbb{Z}) \,$ by:
\begin{equation} \label{a1}
  {\xb} | \g \, := \, \sum_{{\yb} \cdot\check{\g} = {\xb}} \, {\yb}
   \, \in \, \mathbb{Q} [\Sc] \, .
\end{equation}
With this definition one has
\begin{equation} \label{a2}
   \phi_{\xb} | \g \, = \, \phi_{{\xb} | \g} \, ,
   \qquad  \g \in M_2^+ (\mathbb{Z}) \, .
\end{equation}
Modulo the subspace of {\it distribution relations}
$$
 \Rc \, := \, \mathbb{Q}-\text{span of} \quad
 \left\{ {\xb} \, - \, {\xb} | \begin{pmatrix} n &0 \\
 0 &n  \end{pmatrix} \, \, ; \, {\xb} \in {\Sc} \, ,
 n \in \mathbb{Z} \setminus 0 \, \right\} \, ,
$$
the assignment $\, {\xb} \in {\Sc} \longmapsto \phi_{\xb} \,$
induces an isomorphism of $\, \text{PGL}_2^+ (2,
\mathbb{Q})$-modules
\begin{equation} \label{miso}
 \mathbb{Q} [\Sc] / \Rc \, \cong \, \Ec^*_2 (\mathbb Q) \, .
\end{equation}
\medskip

\noindent We are now in a position to state the following result,
which supersedes Lemma \ref{eis}.

\bigskip

\begin{proposition} \label{emu}
For any $\, \g \in M_2^+ (\mathbb{Z})  \,$ one has
\begin{equation}  \label{gb}
 \mu_{\g} \, = \,  2 ( \phi_{\0b} | \g \, - \, \phi_{\0b} ) \, =
 \, 2 \, \left( \sum_{{\yb} \cdot\check{\g} = {\0b}} \, \phi_{\yb}
 \, - \, \phi_{\0b} \right) \, .
\end{equation}
\end{proposition}
\smallskip

\begin{proof}
This follows from Proposition \ref{emu0} and the equalities
 (\ref{a1}), (\ref{a2}).
\end{proof}
\medskip

\noindent The role of the invariant $1$-cocycle needed to complete
the construction of the analogue of the Godbillon-Vey class will
be played by the extension to $ G^+ (\mathbb{Q})$ of the
well-known {\it period cocycle} for modular forms. We start by
recalling its definition.

\noindent Fix $ z_0 \in H $ and, for each $\g \in G^+
(\mathbb{Q})$, define a linear functional $\, \Psi_{z_0} (\g)
\equiv \Psi (\g) \,$ on $ \Om^1  $ by setting
$$
\Psi (\g) \, (\om)  \equiv \langle \om , \Psi (\g) \rangle \, :=
\, \, \int_{z_0}^{\g \cdot z_0} \, \om \, , \qquad \fl \, \om \in
\Om^1 \, .
$$
One easily checks the $1$-cocycle property
$$ \Psi (\g_1 \g_2) (\om) \, = \Psi (\g_2) \, (\g_1^* \om) \, +
\, \Psi (\g_1) (\om) \, , \qquad \fl \, \g_1 , \g_2 \in G^+
(\mathbb{Q}) ,
$$
as well as the fact that the cohomology class $[\Psi] \in  H^{1}
\, ( G^+ (\mathbb{Q}), \, (\Om^1)^* ) \,$ is independent of the
choice of $ z_0 \in H $.

\noindent The cup product of the two $1$-cocycles gives a
$2$-dimensional cohomology class $\, [ {\Eb} \cup \Psi]  \in H^{2}
\, ( G^+ (\mathbb{Q}),
 \, \Om^1 \ot (\Om^1)^*) $. We shall
be interested in the contraction of this class, $ [ \Tr ({\Eb}
\cup \Psi) ]  \in H^{2} \, ( G^+ (\mathbb{Q}), \mathbb{C}) \, ,$
obtained by forming the `trace' 2-cocycle
\begin{equation*}
 \tau (\g_1, \g_2) \equiv
\Tr ({\Eb} \cup \Psi) \, (\g_1, \g_2) \, := \, \langle {\Eb}
(\g_1), \, \g_1 \cdot \Psi (\g_2) \rangle \, , \quad \g_1 , \g_2
\in G^+ (\mathbb{Q}) \, .
\end{equation*}
Explicitly, with $\displaystyle \psi_{\g}:=\frac{\D | \g}{\D}$,
\begin{eqnarray}  \label{tau}
    \tau (\g_1, \g_2) \, &=& \,
\langle \g_1^{-1} \cdot {\Eb} (\g_1) , \, \Psi (\g_2) \rangle \, =
\,
  \langle \g_1^* \left(\  {\Eb} (\g_1) \right),
\, \Psi (\g_2) \rangle  \nonumber \\
\, &=& \, \frac{1}{12  \pi i} \langle \frac{d
\psi_{\g_1}}{\psi_{\g_1}} , \, \Psi (\g_2) \rangle \, = \,
\frac{1}{12  \pi i} \int_{z_0}^{\g_{2}Û \cdot z_0} \, \frac{d
\psi_{\g_1}}{\psi_{\g_1}} \, , \nonumber
\end{eqnarray}
that is
\begin{equation} \label{tau2}
 \tau (\g_1, \g_2) \, = \, \frac{1}{12  \pi i} \cdot
\left( \log \psi_{\g_1} (\g_{2}Û \cdot z_0 ) \, - \, \log
\psi_{\g_1} ( z_0 ) \right) \, ,
\end{equation}
where the determination of the logarithm is unimportant.
\medskip

At this juncture we recall that, as a special case of the main
result in \cite{BY}, the natural map
$$ H_{\rm d}^* (\SL (2, \mathbb{R}),\,  \mathbb{R}) \longrightarrow
H^* (\SL (2, \mathbb{Q}),\,  \mathbb{R})
$$
obtained as the composition of the natural homomorphism
$$
H_{\rm d}^* (\SL (2, \mathbb{R}),\,  \mathbb{R}) \longrightarrow
H^* (\SL (2, \mathbb{R})^{\d},\,  \mathbb{R}) \,
$$
with the restriction map
$$
H^* (\SL (2, \mathbb{R})^{\d},\,  \mathbb{R}) \longrightarrow H^*
(\SL (2, \mathbb{Q}),\,  \mathbb{R})
$$
is an isomorphism. Therefore,
\begin{equation}
H^2 (\SL (2, \mathbb{Q}),\,  \mathbb{R}) \, = \, \mathbb{R} \cdot
e  \, ,
\end{equation}
where $ \, e \in H^2 (\SL (2, \mathbb{Q}),\,  \mathbb{R}) \,$
stands for the Euler class. Noting that the expression (\ref{tau})
actually makes sense for any $\, \g_1, \g_2 \in \GL^+ (2,
\mathbb{R}) \,$, we can view $\tau$ as a $2$-cocycle on $\SL (2,
\mathbb{R})^{\d}$.
\bigskip

\begin{theorem} \label{m1}
    The $2$-cocycle $\, \Re \tau \,
\in \, Z^{2} \, (\SL (2, \mathbb{R})^{\d},\,  \mathbb{R}) \, $
represents the Euler class $\,e \, \in \, H^{2} \, (\SL (2,
\mathbb{R})^{\d}, \, \mathbb{R}) \,$.
\end{theorem}
\medskip

\begin{proof} By construction, one has
\begin{equation*}
\Re \tau (\g_1, \g_2) \, = \, \int_{z_0}^{\g_{2}Û \cdot z_0} \,
\Re (\mu_{\g_1} (z) \, dz) \, , \qquad \fl \, \g_1 , \g_2 \in \
G^+ (\mathbb{R}) := \GL^+ (2, \mathbb{R}) ,
\end{equation*}
Using Proposition \ref{emu0} we may write
\begin{equation} \label{rt}
\Re \tau (\g_1, \g_2) \, = \, \int_{z_0}^{\g_{2}Û \cdot z_0}
  2 \Re (( \phi_{\0b} | \g_1  -  \phi_{\0b} ) (z)  dz)
  \, = \, \int_{z_0}^{\g_{2}Û \cdot z_0} ( \g_{1}^* (\om_0 )  -
   \om_0 ) \, ,
\end{equation}
where
$$  \om_0 (z) \, := \, 2 \Re \left(\phi_{\0b} (z) \, dz\right) \,
= \, \frac{1}{2 \pi^2}
                \, \Re \left( G_{\0b} (z) \, dz \right) \, .
$$
Because $G_{\0b}$ is not holomorphic, $\om_0$ is not a closed
form. Explicitly, by (\ref{GG*}) one has
\begin{equation*} \label{om}
  \om_0 (z) \, = \, \frac{1}{2 \pi^2} \, \Re \left( G_2^* (z)
  \, dz\right) \,
  - \, \frac{1}{2 \pi y}  \, dx \, , \qquad z = x + iy \, ,
\end{equation*}
with $G_2^*$ holomorphic, therefore
\begin{equation} \label{dom}
 d \om_0 (z) \, = \, - \frac{1}{2 \pi} \, \frac{dx \wedge dy}
 {y^2} \, .
\end{equation}
Integrating along geodesic segments in the upper half-plane, we
can rewrite (\ref{rt}) as follows:
\begin{eqnarray} \label{ro}
\Re \tau (\g_1, \g_2)
  \, &=& \, \int_{[z_0 , \, \g_{2}Û \cdot z_0]}  \g_{1}^* (\om_0 ) \, -
   \, \int_{[z_0 , \, \g_{2}Û \cdot z_0]}  \om_0  \, ,
   \nonumber \\
 &=& \, \int_{[\g_{1}Û \cdot z_0 , \, \g_1 \g_{2}Û \cdot z_0]}  \om_0  \, -
   \, \int_{[z_0 , \, \g_{2}Û \cdot z_0]}  \om_0  \, .  \nonumber
\end{eqnarray}
By subtracting the coboundary
\begin{equation*}
 \int_{[z_0 , \, \g_1 \g_{2}Û \cdot z_0]} \om_0  \, -
 \, \int_{[z_0 , \, \g_{1}Û \cdot z_0]}  \om_0  \, -
 \, \int_{[z_0 , \, \g_{2}Û \cdot z_0]}  \om_0
\end{equation*}
we see that $ - \Re \tau$ is cohomologous to the area cocycle
(comp. Appendix B, Remark 4)
\begin{eqnarray}
 \A (\g_1, \g_2)  :&=&
- \int_{[\g_{1}Û \cdot z_0 , \, \g_1 \g_{2}Û \cdot z_0]}  \om_0
\, +
 \int_{[z_0 , \, \g_1 \g_{2}Û \cdot z_0]} \om_0  \, -
 \, \int_{[z_0 , \, \g_{1}Û \cdot z_0]}  \om_0   \nonumber \\
 &=& - \int_{\triangle (z_0 , \g_{1}Û \cdot z_0 , \g_1 \g_{2}Û
 \cdot z_0)}
 d \om_0 \, = \, \frac{1}{2 \pi} \, \text{Area} \left(
 \triangle (z_0 , \, \g_{1}Û \cdot z_0 , \, \g_1 \g_{2}Û \cdot z_0)
 \right) \, . \nonumber
\end{eqnarray}
\end{proof}
\bigskip

\noindent \textbf{Remark 3. \,} One can also show that the
imaginary part $\Im \tau$ is a coboundary.

\noindent Indeed, by (\ref{tau2}) one has, for any $ \g_1 , \g_2
\in \text{PSL} (2, \mathbb{R}) $,
\begin{eqnarray} \label{1t}
12 \pi i \, \tau (\g_1, \g_2) &=& \log \frac{\D | \g_1}{\D} \left(
\g_{2}Û \cdot z_0 \right)
 - \log \frac{\D | \g_1}{\D} \left( z_0 \right)   \nonumber \\
&=& \log \D | \g_1 ( \g_{2}Û \cdot z_0 ) -
\log \D | \g_1 ( z_0 )  \nonumber \\
&-& \left(\log \D ( \g_{2}Û \cdot z_0 ) - \log \D ( z_0 )\right)
\end{eqnarray}
after the cancellation of an additive constant which depends only
on $\g_{1}$.

\noindent Since
$$
\log \D | \g ( z_0 ) \, = \, \log \D (\g \cdot z_0 ) \, - \, 6 \,
\log j^{2} (\g,  z_0 ) \, + \, 2 \pi i \, k (\g) \, ,
$$
for some $\,  k (\g) \in \mathbb{Z} \,$ and with the choice of the
principal branch for $\, \mathbb{C} \setminus [0, \infty) \,$ of
the logarithm of the square of the automorphy factor, the
succession of identities in (\ref{1t}) can be continued with
\begin{eqnarray} \label{2t}
12 \pi i \, \tau (\g_1, \g_2) &=&
 \log \D (\g_ 1 \g_2 \cdot z_0 ) - \log \D (\g_1 \cdot z_0 )
 - \log \D (\g_2 \cdot z_0 )  \nonumber \\
 &+& \log \D ( z_0 )
\, -\,  6 \left( \log j^{2} (\g_1,  \g_2 \cdot z_0 ) - \log j^{2}
(\g_1, z_0 ) \right) \, .
\end{eqnarray}
The equality
\begin{equation} \label{log}
 \log j^{2} (\g_1 \g_2 ,  z_0 ) =  \log j^{2} (\g_1,  \g_2 \cdot z_0 )
 +  \log j^{2} (\g_2 ,  z_0 ) + 2 \pi i \, c (\g_1,  \g_2 ) \, ,
\end{equation}
determines a cocycle $c \in Z^2 (\PSL (2, \mathbb{R}), \,
\mathbb{Z}) \,$, which is precisely the cocycle discussed
in~\cite[\S B-2]{BG}. Inserting (\ref{log}) into (\ref{2t}) one
obtains
\begin{eqnarray} \label{3t}
12 \pi i \, \tau (\g_1, \g_2) &=&  \log \D ( z_0 )     \nonumber \\
&+& \log \D (\g_ 1 \g_2 \cdot z_0 ) - \log \D (\g_1 \cdot z_0 )
 - \log \D (\g_2 \cdot z_0 ) \nonumber \\
&-& 6 \left( \log j^{2} (\g_1 \g_2 \cdot z_0 ) - \log j^{2} (\g_2,
z_0 )
- \log j^{2} (\g_1, z_0 ) \right) \nonumber \\
&+& 12 \pi i \, c (\g_1,  \g_2 ) \, .
\end{eqnarray}
This shows that the cocycles $\, \tau \,$ and $\, c \,$ differ by
a coboundary, which proves our claim. At the same time, in
conjunction with the above theorem, it gives an alternate proof to
the statement (cf.~\cite[Lemma 2.1]{BG}) that the cocycle $ \, c
\,$ represents the Euler class.

\bigskip

We shall now refine the above construction to produce a remarkable
rational representative for the Euler class $\,e \, \in \, H^{2}
\, (\SL (2, \mathbb{Q}), \, \mathbb{Q}) $. In the process, we
shall make extensive use of the results in \cite[Chap. 2]{St}. To
begin with, we replace the $1$-cocycle $\, \Psi = \Psi_{z_0} \,$
by a cohomologous $1$-cocycle -- introduced by Stevens in
\cite[Def. 2.3.1]{St} -- which is independent of the base point
$\, z_0 \in H \, $. The new cocycle, $\Pi \in  Z^{1} \, ( G^+
(\mathbb{Q}), \, (\Mc_2)^* ) \,$, is given by the formula
$$
\Pi (\g) \, (\om) \, := \, \, \int_{z_0}^{\g \cdot z_0} f(z) dz
\, - \, z_0 \cdot
  {\ab}_0 (f|\g - f) \, + \,
\int_{z_0}^{i \ify} \left( \widetilde{f|\g} - \widetilde{f}
 \right) (z) \, dz \, ,
$$
where  $ \, f \in \Mc_2 (\G (N)) \, $ for some level $N$,
$$ f(z) \, = \, \sum_{n=0}^{\ify} \, {\ab}_n  \,e^{\frac{2\pi i n z}{N}}
$$
represents its Fourier expansion at $\ify$, and
$$
 \widetilde{f} (z) \, := \, f(z) \, - \, \ {\ab}_0 (f) \, \, ;
$$
note that $\ {\ab}_0 (f) \,$, and therefore $\, \widetilde{f} \,$
too, are independent of the level.

\noindent By coupling ${\Eb}$ with  $\Pi$ instead of $\Psi$, one
gets a 2-cocycle equivalent to the restriction of $\tau $ to $G^+
(\mathbb{Q})$,
$$ \theta (\g_1, \g_2) :=
  \langle {\Eb} (\g_1), \, \g_1 \cdot
\Pi (\g_2) \rangle \, , \qquad \g_1 , \g_2 \in G^+ (\mathbb{Q}) \,
\, ;
$$
explicitly,
\begin{eqnarray} \label{theta}
   \theta (\g_1, \g_2)  &=&
\int_{z_0}^{\g_{2}Û \cdot z_0}\mu_{\g_1} (z) dz \, - \, z_0 \cdot
 {\ab}_0 (\mu_{\g_1}|\g_2 - \mu_{\g_1} ) \nonumber \\
&+&\, \int_{z_0}^{i \ify} ( \widetilde{\mu_{\g_1}|\g_2} -
\widetilde{\mu_{\g_1}}
 ) (z) \, dz  .
\end{eqnarray}
Furthermore, as follows from the preceding proof, cf. (\ref{3t}),
retaining only the real part of $\, \theta \,$ gives a
cohomologous cocycle
\begin{equation} \label{rho}
\rho (\g_1, \g_2) := \Re \, \theta (\g_1, \g_2)  \, = \,
 \Re \, \langle {\Eb} (\g_1), \, \g_1 \cdot
\Pi (\g_2) \rangle , \quad \g_1 , \g_2 \in G^+ (\mathbb{Q}) \, .
\end{equation}

\noindent A priori, $ \rho \in  Z^{2} \, ( G^+ (\mathbb{Q}),
\mathbb{R} ) $, and we want to show that it actually takes
rational values, which moreover can be explicitly computed. To
this end, we now proceed to unfold the expression of the cocycle
$\rho (\g_1, \g_2) = \Re \theta (\g_1, \g_2) $, defined by
(\ref{rho}), for
\begin{equation} \label{gg}    \g_1 \, = \, \begin{pmatrix} a_1 &b_1 \\
 c_1 &d_1  \end{pmatrix} \, , \quad
\g_2 \, = \, \begin{pmatrix} a_2 &b_2\\
 c_2 &d_2  \end{pmatrix} \,  \in \, M_2^+ (\mathbb{Z})  \, .
\end{equation}

\noindent The expression of $\, \rho (\g_1, \g_2) \, $ is simpler
when
$$
 \g_2 \, \in \, B^+ (\mathbb{Z}) \, , \quad \text{i.e.} \quad
 c_2 \, = \, 0 \, .
$$
Indeed, in that case if we let $\, z_0 \ra i \ify \, $  then $\,
\g_2 \cdot z_0 \ra i \ify \, $ as well, hence both integrals in
the right hand side of (\ref{theta}) will tend to $0$. On the
other hand, the limit of the middle term can be easily computed,
cf. \cite[\S 2.3]{St}, giving
$$
\theta (\g_1, \g_2)  \, = \,
 \frac{b_2}{d_2} \, {\ab}_0 (\mu_{\g_1} )    \, .
$$
Thus, by (\ref{gb})
$$ \rho (\g_1, \g_2) \, = \,  \frac{b_2}{d_2} \cdot
  2 \, {\ab}_0  ( \phi_{\0b} | \g_1 \, - \, \phi_{\0b} ) \,    \, .
$$
Since according to \cite[Prop. 2.5.1 (ii)]{St},
\begin{equation} \label{ct}
{\ab}_0  ( \phi_{\xb} ) \, = \, \frac{1}{2} \, {\Bb}_{2} (x_1) \,
,
\end{equation}
from (\ref{a1}), (\ref{a2}) it follows that, when $\, c_2 = 0 $,
\begin{equation} \label{f1}
  \rho (\g_1, \g_2) \, = \,  \frac{b_2}{d_2} \,
 \left( \sum_{{\yb} \cdot \check{\g_1} = {\0b}} \,  {\Bb}_{2} (y_1)
 - \, \frac{1}{6} \right) \, .
\end{equation}
Here $ \, {\Bb}_2 \, : \, \mathbb{Q}\slash \mathbb{Z} \, \ra \,
\mathbb{Q} \, $ is the periodized function corresponding to the
Bernoulli polynomial $ \, B_2 (X) \, = \, X^2  - X + \frac{1}{6}
\, $ so that $ \, {\Bb}_2 (x):= B_2 (x-[x])$ where $[x]$
is the integer part of $x$. 
\medskip

\noindent The complementary case, when
$$
\g_2 \, = \, \begin{pmatrix} a_2 &b_2\\
 c_2 &d_2  \end{pmatrix} \,  \in \, M_2^+ (\mathbb{Z}) \,
 \quad \hbox{with} \quad
 c_2 \, > \, 0 \, ,
$$
can be dealt with by means of Stevens' {\it Dedekind symbol} map
\cite[Def. 2.5.2]{St} $\, {\Sb}_E : \mathbb{P}^1 (\mathbb{Q}) \ra
\mathbb{Q} \,$ associated to an Eisenstein series $ E \in \Ec_2
(\mathbb{Q}) $. According to \cite[Prop 2.5.4]{St}, the Dedekind
symbol of a basis element $\, \phi_{\xb} \,$ is given by the
following generalized Dedekind sum, introduced by Meyers in
\cite{Me}: for $ \quad m, n \in \mathbb{Z} \quad \hbox{with}
 \quad (m, n) = 1 \quad \hbox{and} \quad n > 0 \,$,
\begin{equation} \label{ded}
 {\Sb}_{\phi_{\xb}} (\frac{m}{n}) \, = \, \sum_{j=0}^{n-1}
 {\Bb}_1 \left(\frac{x_1 + j}{n} \right) \,
 {\Bb}_1 \left(\frac{m (x_1 + j)}{n} + x_2 \right) \, ,
\end{equation}
where $ \, {\Bb}_1 : \, \mathbb{Q}\slash \mathbb{Z} \, \ra \,
\mathbb{Q} \, $ denotes the periodized function corresponding to
the Bernoulli polynomial $ \, B_1 (X) \, = \, X - \frac{1}{2}  \,
$. Relying again on the results in \cite[\S 2.5]{St}, in
particular on formula (2.5.3), one obtains:
\begin{equation*}
\rho (\g_1, \g_2) \, = \, \frac{a_2}{c_2} \, {\ab}_0 (\mu_{\g_1} )
\, + \, \frac{d_2}{c_2} {\ab}_0 (\mu_{\g_1} | \g_2 ) \, - \,
{\Sb}_{\mu_{\g_1}} ( \frac{a'_2}{c'_2}) \, ,
\end{equation*}
where $ \,  \frac{a'_2}{c'_2} \, , \, (a'_2 , c'_2) = 1 \,$, is
the fraction $ \,  \frac{a_2}{c_2} \,$ in lowest terms form. Using
once more (\ref{gb}), as well as (\ref{ct}) and (\ref{ded}), one
finally obtains, for the case $ \, c_2 > 0 $,
\begin{eqnarray} \label{f2}
\rho (\g_1, \g_2)  &=& \frac{a_2}{c_2} \, \left( \sum_{{\yb} \cdot
\check{\g_1} = {\0b}}  {\Bb}_{2} (y_1)
 - \, \frac{1}{6} \right)  \nonumber \\
&+& \frac{d_2}{c_2} \left( \sum_{{\zb} \cdot
 \check{\g_2} \check{\g_1} = {\0b}}
{\Bb}_{2} (z_1) \, - \, \sum_{{\yb} \cdot \check{\g_2} = {\0b}}
{\Bb}_{2} (y_1) \right)   \nonumber \\
&-& 2 \sum_{{\yb} \cdot \check{\g_1} = {\0b}} \sum_{j=0}^{c'_2 -1}
 {\Bb}_1 \left(\frac{y_1 + j}{c'_2} \right)
 {\Bb}_1 \left(\frac{a'_2 (y_1 + j)}{c'_2} + y_2 \right) \nonumber \\
&+& 2 \sum_{j=0}^{c'_2 -1}
 {\Bb}_1 \left(\frac{j}{c'_2} \right)
 {\Bb}_1 \left(\frac{a'_2 \, j}{c'_2} \right) \, .
 \end{eqnarray}

\noindent In conclusion, we have proved:
\bigskip

\begin{theorem} \label{m2}
    The $2$-cocycle $\, \rho \,
\in \, Z^{2} \, (\PSL (2, \mathbb{Q}),\,  \mathbb{Q}) \, $
given by the formulae (\ref{f1}) and (\ref{f2}) represents the
Euler class $\,e \, \in \, H^{2} \, (\PSL (2,
\mathbb{Q}),\, \mathbb{Q})$.
\end{theorem}
\bigskip

\section*{Appendix A -- Hopf cyclic cohomology}

\noindent We shall recall in this appendix the definition of
cyclic cohomology for Hopf algebras endowed with a modular pair in
involution ($\nu, \sigma$), in the special case $\sigma=1$ (cf.
\cite{CM1}. We refer to \cite{CM2} for the general case).
\medskip

\noindent Let $\Hc$ be a Hopf algebra, $\nu$ a character of $\Hc$
such that the twisted antipode
 $\wt S = \nu * S$ satisfies the involutive property
\begin{equation} \nonumber
\wt S^2 = \Id \, .
\end{equation}
 The cyclic cohomology groups
$H  C_{\Hopf}^* (\Hc)$ are defined by means of the cyclic module
associated to the Hopf algebra
 $\Hc$ as follows. For each $\, n \in \mathbb{N} \,$,
$ \, C^n (\Hc) = \Hc^{\ot n} \, \, ;$ the face operators
$\partial_i: C^{n-1} (\Hc) \ra C^n (\Hc), \quad 0 \leq i \leq n \,
,$ are
\begin{eqnarray}
\partial_0 (h^1 \ot \ldots \ot h^{n-1}) \, &= \, 1 \ot h^1
\ot \ldots \ot h^{n-1}, \qquad \qquad \qquad
\qquad \qquad \nonumber \\
\partial_j (h^1 \ot \ldots \ot h^{n-1}) \, &= \, h^1 \ot \ldots \ot \D h^j \ot
\ldots \ot h^{n-1} \, , \quad 1 \leq j \leq n-1 \, , \nonumber \\
\partial_n (h^1 \ot \ldots \ot h^{n-1}) \, &= \, h^1 \ot \ldots \ot h^{n-1}
\ot 1 \, \, ;  \qquad \qquad \qquad \qquad \nonumber
\end{eqnarray}
the degeneracy operators $\s_i : C^{n+1} (\Hc) \ra C^n (\Hc), \,
\quad  0 \leq i \leq n \, ,$ are
\begin{eqnarray}
\s_i (h^1 \ot \ldots \ot h^{n+1}) \, = \, h^1 \ot \ldots \ot \ve
(h^{i+1}) \ot \ldots \ot h^{n+1}  \, \, ; \nonumber
\end{eqnarray}
finally, the cyclic operator $\, \tau_n : C^n (\Hc) \ra C^n (\Hc)
\,$ is given by
\begin{eqnarray}
    \tau_n (h^1 \ot \ldots \ot h^n) \, = \, (\D^{n-1}  \wt S (h^1)) \cdot h^2
\ot \ldots \ot h^n \ot 1   \, . \nonumber
\end{eqnarray}
The Hopf cyclic cohomology is computed from the normalized
bicomplex $\, (C C^{*, *} (\Hc), \, b , \, B )$, where:
\begin{eqnarray*} \label{CbB}
CC^{p, q} (\Hc) &=& {\bar C}^{q-p} (\Hc), \quad q \geq p , \cr
 C C^{p, q} (\Hc) &=& 0 ,  \quad q < p \, ; \cr
\end{eqnarray*}
with
$$ {\bar C}^n (\Hc) =
 \cap \;{\rm Ker}\,\s_i \,, \quad \fl n \geq 1, \quad \quad
                        {\bar C}^0 (\Hc) = \mathbb{C} ;$$
 The operator
\begin{equation*}  \nonumber
b: {\bar C}^{n-1} (\Hc) \ra {\bar C}^n (\Hc), \qquad b =
\sum_{i=0}^{n} (-1)^i \partial_i \,
\end{equation*}
has the expression
\begin{eqnarray} \nonumber
\ b (h^1 \ot \ldots \ot h^{n-1}) &=& 1 \ot h^1 \ot \ldots \ot
h^{n-1}  \cr & + & \ \sum_{j=1}^{n-1} (-1)^j \sum_{(h_{j})} h^1
\ot \ldots \ot h_{(1)}^j \ot h_{(2)}^j \ot \dots \ot h^{n-1} \cr
&+ & \ (-1)^n h^1 \ot \ldots \ot h^{n-1} \ot 1  ,
\end{eqnarray}
while for $n=0$ , $\, b( \mathbb{C} ) \, = \, 0 \, $.

\noindent The $B$-operator $B: {\bar C}^{n+1} (\Hc) \ra {\bar C}^n
(\Hc)$ is defined by the formula
\begin{equation} \nonumber
  B = A \circ B_{0} \, ,  \quad n \geq 0 \, ,
\end{equation}
where $B_{0} : {\bar C}^{n+1} (\Hc) \ra {\bar C}^n (\Hc) \,$ is
the (extra degeneracy) operator
\begin{eqnarray} \nonumber
&& B_{0} (h^1 \ot \ldots \ot h^{n+1}) = (\D^{n-1}
  \wt S (h^1)) \cdot h^2
\ot \ldots \ot h^{n+1} \nonumber \\
&& \quad = \sum_{(h^{1})} S(h^1_{(n)}) h^2 \ot \ldots \ot
S(h^1_{(2)}) h^n \ot \wt S (h^1_{(1)})  h^{n+1} , \nonumber
\end{eqnarray}
$$ B_{0} (h) = \nu (h), \, \, h \in \Hc $$
and
\begin{equation} \nonumber
A = 1 + (-1)^n \tau_n  + \ldots + (-1)^{n^{2}} {\tau_n}^n \, .
\end{equation}
The groups $\, H  C_{\Hopf}^n (\Hc) \,$ are computed from the
first quadrant total complex $\, ( TC^{*} (\Hc),\, b+B ) \,$,
\begin{equation*}
    TC^{n}(\Hc) = \sum_{p=0}^{n} \, C C^{p, n-p} (\Hc) \, ,
\end{equation*}
while the ($\mathbb{Z} /2$--graded) periodic groups $\, PH
C_{\Hopf}^* (\Hc) \,$ are computed from the full total complex $\,
(P T C^{*}(\Hc),\,  b+B ) \,$,
\begin{equation*}
      P T C^{n}(\Hc) = \sum_{p} \, C C^{p, n-p} (\Hc) \, .
\end{equation*}
\bigskip

\section*{Appendix B -- Godbillon-Vey as Hopf cyclic class}

After recalling its original definition, we shall give in this
appendix the promised detailed account of the interpretation of
the Godbillon-Vey as Hopf cyclic class for $\Hc_{1}$.
\medskip

\noindent Let $V$ be a closed, smooth manifold, foliated by a
transversely oriented codimension $1$ foliation $\Fc$. Then $T\Fc
= \Ker \om \sbs TV$, for some $\om \in \Om^1 (V)$ such that $\om
\wdg d\om = 0$. Equivalently, $d\om = \om \wdg \a$ for some $\a
\in \Om^1 (V)$, which implies $d\a \wdg \om = 0$. In turn, the
latter ensures that $d\a = \om \wdg \b$, $\b \in \Om^1 (V)$. Thus,
$\a \wdg d\a \in \Om^3 (V)$ is closed. Its de~Rham cohomology
class,
\begin{equation*}
\GV (V, \Fc) = [\a \wdg d\a] \in H^3 (V,\mathbb{R}) \, ,
\end{equation*}
is independent of the choices of $\om$ and $\a$ and represents the
original definition of the Godbillon-Vey class.
\smallskip

\noindent The Godbillon-Vey class acquires a universal status when
viewed as a characteristic class (cf. \cite{Ha}) associated to the
Gelfand-Fuchs cohomology \cite{GF} of the Lie algebra ${\Fa}_1 =
\mathbb{R} \, [[x]] \, \partial_x$ of formal vector fields on
$\mathbb{R}$. The cohomology $\, H^{*} ({\Fa}_1 , \mathbb{R}) \, $
is finite dimensional and the only nontrivial groups are:
$$
H^0 ({\Fa}_1 , \mathbb{R}) = \mathbb{R} \cdot 1 \qquad \hbox{and}
\qquad H^3 ({\Fa}_1 , \mathbb{R}) = \mathbb{R} \cdot gv \, ,
$$
where
\begin{equation*}
gv (p_1 \, \partial_x , p_2 \, \partial_x , p_3 \, \partial_{x}) =
\left\vert \begin{array}{ccc}
p_1 (0) &p_2 (0) &p_3 (0) \\ \\
p'_1 (0) &p'_2 (0) &p'_3 (0) \\ \\
p''_1 (0) &p''_2 (0) &p''_3 (0)
\end{array} \right\vert \, ,
\end{equation*}
that is, with the obvious notation,
\begin{equation*}
gv = \t^0 \wdg \t^1 \wdg \t^2 \, .
\end{equation*}

\noindent Given any oriented $1$-dimensional manifold $M^1$, the
Lie algebra cocycle $gv$ can be converted into a $3$-form on the
jet bundle (of orientation preserving jets)
$$ J_+^{\ify} (M^1) = \limproj_n J_+^n (M^1) \, ,
$$
invariant under the pseudogroup $\Gc^+ (M^1)$ of all orientation
preserving local diffeomorphisms of $M^1$. Indeed, sending the
formal vector field
$$
p = j_0^{\ify} \left( \frac{dh_t}{dt} \Biggl\vert_{t=0} \right)
\in \Fa_1 \, ,
$$
where $\{ h_t \}$ is a $1$-parameter family of local
diffeomorphisms of $\mathbb{R}$ preserving the origin, to the
$\Gc^+ (M^1)$-invariant vector field
$$
j_0^{\ify} \left( \frac{d(f \circ h_t)}{dt} \Biggl\vert_{t=0}
\right) \in T_{j_0^{\ify} (f)} \, J_+^{\ify} (M^1) \,
$$
gives a natural identification of the Lie algebra complex of
${\Fa}_1$ with the invariant forms on the jet bundle,
$$
 \quad \t \in C^{\bu} ({\Fa}_1) \mpo \t \in \Om^{\bu}
(J_+^{\ify} (M^1))^{\Gc^+ (M^1)} \, .
$$
In local coordinates on $J_+^{\ify} (M^1)$, given by the
coefficients of the Taylor expansion at $0$,
$$
f(s) = y + s \, y_1 + s^2 y_2 + \cdots \, , \qquad y_1 > 0 \, ,
$$
one has
\begin{eqnarray}
dy \, &=& \, y_1 \, \t^0  \nonumber \\
dy_1 \, &=& \, y_1 \, \t^1 + 2 \, y_2 \, \t^0
                        \nonumber \\
dy_2 \, &=& \, y_1 \, \t^2 + 2 \, y_2 \, \t^1 + 3 \, y_3 \, \t^0
\, , \nonumber
\end{eqnarray}
therefore
\begin{equation*}
 \quad gv = \frac{1}{y_1^3} \, dy \wdg dy_1 \wdg dy_2 \in
\Om^3 (J_+^{\ify} (M^1))^{\Gc^+ (M^1)} \, .
\end{equation*}
\medskip

\noindent Given a codimension $1$ foliation $(V, \Fc)$ as above,
one can find an open covering $\{ U_i \}$ of $V$ and submersions
$f_i : U_i \ra T_i \sbs \mathbb{R}$, whose fibers are plaques of
$\Fc$, such that $f_i = g_{ij} \circ f_j$ on $U_i \cap U_j$, with
$g_{ij} \in \Gc^+ (M^1)$ a $1$-cocycle. Then $\, M^1 = \bigcup_i
T_i \times \{i\} \, $ is a complete transversal. Let $J^{\infty}
(\Fc)$ denote the bundle over $V$ whose fiber at $x \in U_i$
consists of the $\infty$-jets of local submersions of the form
$\varphi \circ f_i$ with $\varphi \in \Gc^+ (M^1)$. Using the
$\Gc^+ (M^1)$-invariance of $gv \in \Om^3 (J^{\infty} (M^1))$ one
can pull it back to a closed form $gv (\Fc) \in \Om^3 (J^2
(\Fc))$. Its de~Rham class  $ [gv (\Fc)] \in H^3 (J^2 (\Fc),
\mathbb{R})$, when viewed as a class in $H^3 (V, \mathbb{R})$ (the
fibers of $J^2 \Fc$ being contractible), is precisely the
Godbillon-Vey class $\GV (V, \Fc)$.

In \cite{CM1, CM2} we proved in full generality that the Hopf
cyclic cohomology of the Hopf algebra of transverse differential
operators on an $n$-dimensional manifold is canonically isomorphic
to the Gelfand-Fuchs cohomology of the Lie algebra of formal
vector fields on $\mathbb{R}^{n}$, via an isomorphism which can be
explicitly realized at the cochain level in terms of a fixed
torsion-free connection. In particular, one has a canonical
isomorphism
\begin{equation*}
\kappa_{1}^{*} \, : \, H^{*} ({\Fa}_1 , \mathbb{C})
\build\longrightarrow_{}^{\simeq}P HC_{\Hopf}^{*} \, (\Hc) \, .
\end{equation*}
\smallskip

\begin{proposition} \label{gv1} The canonical cochain map
associated to the trivial connection on $\, J_+^1 (\mathbb{R}) \,$
sends the universal Godbillon-Vey cocycle $\, gv \,$  to the Hopf
cyclic cocycle $ \, \d_{1} \,$, implementing the identity
\begin{equation*}
\kappa_{1}^{*} ([gv]) \, = \, [\d_{1}] \, .
\end{equation*}
\end{proposition}
\smallskip

\begin{proof} To begin with
recall that as a form on the jet bundle $J_+^2 (\mathbb{R})$
$$
gv \, = \, \frac{1}{y_1^3} \, dy \wdg dy_1 \wdg dy_2 \, .
$$
According to the two-step definition in \cite{CM1}, \cite{CM3}  of
the isomorphism $\kappa_{1}^{*}$, one first turns the Lie algebra
cocycle $gv \in C^{2} ({\Fa}_1, \mathbb{R})$ into a group
$1$-cocycle $\, C_{1,0} (gv) \, $ on $\, \Gc = \Diff^+
(\mathbb{R}) \,$ with values in currents on $\, J_+^1 (\mathbb{R})
\,$, and then one takes its image in the cyclic bicomplex 
 under the canonical map $\Phi$. The resulting cyclic
cocycle
$$
(\Phi \, ( C_{1,0} \, (gv) \,) ) \, (f^0 \, U_{\vp_0}^* , f^1 \,
U_{\vp_1}^*)
$$
is automatically supported at the identity, i.e. it is nonzero
only when $\vp_1 \vp_0 = 1$. Moreover, it is of the form
$$
(\Phi \, ( C_{1,0} \, (gv) \, ) ) \, (f^0 \, U_{\vp}^* , f^1 \,
U_{\vp^{-1}}^*) \, = \, -
  \langle C_{1,0} (gv) (1, \vp) , \, f^0 \cdot f^1 \circ \vp \rangle \, .
$$
By definition,
\begin{equation*}
\langle C_{1,0} (gv) (1 , \vp) , f \rangle \, = \, \int_{\D^1 \ts
J_+^1 (\mathbb{R})} f \, \wt{\s} (1 , \vp)^* \left( gv \right)
\end{equation*}
where $\, \D^1 \,$ is the $1$-simplex and
$$ \wt{\s} (1 , \vp)  : \D^1 \ts J_+^1 (\mathbb{R} )
  \ra J_+^{\ify} (\mathbb{R})
$$
has the expression
$$
\wt{\s} (1 , \vp) (t , y , y_1) \, = \, \s_{(1-t) \nb_{0} + t
\nb_{0}^{\vp}} (y , y_1) \, ,
$$
whose meaning we now proceed to explain.

First, $\, \nb_{0} \,$ stands for the trivial linear connection on
$\mathbb{R}$, given by the connection form on  $\, J_+^1
(\mathbb{R}) \,$
$$
\om_{0} = y_1^{-1} \, dy_1 \, ,
$$
while $\, \nb_{0}^{\vp} \,$  denotes its transform under the
prolongation
$$
\vp (y, y_1) = (\vp (y) , \, \vp' (y) \cdot y_1) \, ,
$$
of the diffeomorphism $ \vp \in  \Gc $\, ; the latter corresponds
to the connection form
$$
\vp^* (\om_{0}) = \frac{1}{\vp' (y) \, y_1} \, (\vp' (y) \, dy_1 +
\vp'' (y) \cdot y_1 \, dy) = y_1^{-1} \, dy_1 + \frac{d}{dy} (
\log \vp' (y)) \, dy \, .
$$
Furthermore, for any linear connection $\, \nb \,$ on
$\mathbb{R}$, $ \, \s_{\nb} \,$ denotes the jet
$$
\s_{\nb} (y, y_1) = j_0^{\ify} \left( Y (s) \right) \,
$$
of the local diffeomorphism
$$ s \, \mpo  \, Y(s) := \exp_y^{\nb} \left( s \, y_1 \, \frac{d}{dy}
\right) \, , \quad s \in \mathbb{R} \, .
$$
Now $ \, Y(s) \,$ satisfies the geodesics ODE
$$
\left\{ \begin{matrix} \ddot Y (s) + \G_{1 \,1}^{1} (Y(s)) \cdot
\dot Y (s)^2 = 0 \, , \hfill \cr Y(0) = y \, , \hfill \cr \dot Y
(0) = y_1 \, . \hfill \cr
\end{matrix} \right. \,
$$
Since we only need the 2-jet of the exponential map, suffices to
retain that
$$
Y (0) = y \, , \quad \dot Y (0) = y_1 \quad \hbox{and} \quad \ddot
Y (0) = - \G_{1 \,1}^{1} (y) \, y_1^2 \, .
$$
Thus,
$$
\s_{\nb} (y, y_1) \, = \, y + y_1 \, s \, - \, \G_{1 \,1}^{1} (y)
\, y_1^2 \, s^2  \, + \, \hbox{higher order terms} \, .
$$

In our case  $\, \nb \, = \, (1-t) \nb_{0} \,  + \, t
\nb_{0}^{\vp} \, $, which gives
$$
\G_{1 \,1}^{1} (t, y) \, = \, t \, \frac{d}{dy} ( \log \vp' (y))
\, ,
$$
and therefore
$$
\wt{\s} (1 , \vp) (t , y , y_1) \, = \, \, y + y_1 \, s \, - \, t
\, \frac{d}{dy} ( \log \vp' (y))
 \, y_1^2 \, s^2  \, + \, \hbox{higher order terms} \, .
$$
It follows that
$$
\wt{\s} (1 , \vp)^* (dy) \, = \, dy \, , \qquad \wt{\s} (1 ,
\vp)^* (dy_1) \, = \, dy_1
$$
and
\begin{eqnarray}
\wt{\s} (1 , \vp)^* (dy_2) \, = \, - \left(  \frac{d}{dy} ( \log
\vp' (y)) \, dt \, + \, t \frac{d}{d y} \left(\frac{d}{dy} ( \log
\vp' (y))\right)
\, dy \right) y_1^2  \nonumber \\
 - \, 2 t \, \frac{d}{dy} ( \log \vp' (y)) \, y_1
\, d y_1 \, . \nonumber
\end{eqnarray}
Hence on $\D^1 \ts J_+^1 (\mathbb{R})$,
$$
\wt{\s} (1 , \vp)^*  (gv) \,  = \, - \frac{1}{y_1} \, \frac{d}{dy}
 ( \log \vp' (y))
\, dt \wdg dy \wdg dy_1 \, .
$$
Going back to the definition of the group cochain, one gets
\begin{eqnarray}
\langle C_{1,0} (gv) (1 , \vp) , f \rangle \, &=& - \int_{J_+^1
(\mathbb{R})} f (y , y_1) \cdot \int_0^1
 dt \cdot \frac{1}{y_1} \, \frac{d}{dy} ( \log \vp' (y))
 dy \wedge dy_1 \nonumber
\\
&= & - \int_{J_+^1 (\mathbb{R})} f(y , y_1) \left( y_{1} \
\frac{d}{dy} ( \log \vp' (y)) \right) \, \frac{dy \wedge
dy_1}{y_1^{2}} \nonumber \, ,
\end{eqnarray}
which finally gives
\begin{eqnarray}
(\, \Phi \, ( C_{1,0} \, (gv))&)&
(f^0 \, U_{\vp}^* , f^1 \, U_{\vp^{-1}}^*) \, = \nonumber \\
&=& \, \int_{J_+^1 (\mathbb{R})} f^0 \cdot f^1 \circ \vp \cdot
\left( y_1 \, \frac{d}{dy} ( \log \vp' (y)) \right) \,
\frac{dy \wedge dy_1}{y_1^{2}} \nonumber \\
= \tau (f^0 U_{\vp}^* \cdot \d_1 (f^1 U_{\vp^{-1}}^*)) &=& \,
\chi_{\tau} (\d_1) (f^0 \, U_{\vp}^* , f^1 \, U_{\vp^{-1}}^*) \, .
\nonumber
\end{eqnarray}
\end{proof}
\bigskip

\noindent {\bf Remark 4. \,} Since
\begin{eqnarray}
 C_{1,0} (gv) (1 , \vp) \, &=& \,
 - \left( y_{1} \
\frac{d}{dy} ( \log \vp' (y)) \right) \cdot
\frac{dy \wedge dy_1}{y_1^{2}}  \nonumber \\
&=&\, - \frac{d}{dy} ( \log \vp' (y)) \frac{dy \wedge dy_1}{y_1}
\, = \, d \left(\frac{d}{dy} (\log \vp' (y)) \, \log y_1 \, dy
\right) \, , \nonumber
\end{eqnarray}
it follows (as in ~\cite[ pp. 41-42]{Ha2}) that, in the
(non-homogeneous) $(d, \d)$-bicomplex
$$\{ \, \G^{*,*} \left( \Diff^+ (S^1), J_+^1 (S^1) \right) , \,
d, \, \d \, \}
$$
computing the $\Diff^+ (S^1)$-equivariant cohomology of $J_+^1
(S^1)$,
$$ \mu (\vp) \, := \, C_{1,0} (gv) (1 , \vp)
$$
is cohomologous to the $\d$-boundary of
$$
  \nu (\vp) \, = \, \frac{d}{dy} (\log \vp' (y)) \, \log y_1 \, dy \, ,
$$
that is with
\begin{equation*}
\d \nu \, (\vp, \psi) \, = \,
 d (\log \vp' \circ \psi) \cdot\log \psi' \, .
\end{equation*}
The latter is precisely Thurston's formula for the Godbillon-Vey
class. By integration along the fiber one gets the Bott-Thurston
cocycle (see~\cite{Bo}) on $\Diff(S^{1})$
\begin{equation*}
\BT \, (\vp , \psi) \, = \, \int_{S^{1}} \log \psi' \cdot d \log
(\vp' \circ \psi) \, .
\end{equation*}
In turn, its restriction to $\text{PSL}(2, \mathbb{R})$, acting by
fractional linear transformation on the projective real line
$\mathbb{P}^1(\mathbb{R})$, gives a multiple of the \textit{area
cocycle} representing the Euler class (cf.~\cite[{\it Appendix} by
R. Brooks]{Bo}):
\begin{equation*}
  \A (\g_1 , \g_2) \, = \, \text{Area}
 \left(\triangle(i, \, \g_2 \cdot i, \, \g_2 \g_1 \cdot i)\right) \, ,
 \qquad  \g_1 , \g_2 \in  \PSL(2, \mathbb{R}) \, ,
\end{equation*}
where $\text{Area} \left(\triangle(z_1, z_2, z_3 )\right)$ denotes
the hyperbolic area of the geodesic triangle in the upper
half-plane with vertices $z_1, z_2, z_3$.

\end{document}